\documentclass[amsbook,12pt]{article}
\usepackage{epsfig, amsfonts,amssymb, amsmath}
\usepackage[usenames]{color}
\usepackage{picinpar}
\usepackage{float, graphicx}
\newcommand{\nl}{\mbox{}\\}

\begin{document}
%
%
%
\mbox{} \vspace{-2.000cm} \\
\begin{center}
\mbox{\Large \bf %
Some basic properties of bounded solutions} \\
\mbox{} \vspace{-0.250cm} \\
{\Large \bf %
of parabolic equations with \mbox{\boldmath $p$}-{\large L}aplacian diffusion} \\
\nl
\mbox{} \vspace{-0.350cm} \\
{\large \sc J.\;Q.\;Chagas,}$\mbox{}^{\!\:\!1}$
{\large \sc P.\;L.\;Guidolin}$\mbox{}^{\;\!2}$
{\large \sc and J.\;P.\;Zingano}$\mbox{}^{\;\!3}$ \\
\mbox{} \vspace{-0.125cm} \\
$\mbox{}^{1}${\small
Departamento de Matem\'atica e Estat\'\i stica} \\
\mbox{} \vspace{-0.685cm} \\
{\small
Universidade Estadual de Ponta Grossa} \\
\mbox{} \vspace{-0.685cm} \\
{\small
Ponta Grossa, PR 84030-900, Brazil} \\
\mbox{} \vspace{-0.350cm} \\
$\mbox{}^{2}${\small
Instituto Federal de Educa\c c\~ao, Ci\^encia e Tecnologia} \\
\mbox{} \vspace{-0.685cm} \\
{\small
Farroupilha, RS 95180-000, Brazil} \\
\mbox{} \vspace{-0.350cm} \\
$\mbox{}^{3}${\small Departamento de Matem\'atica Pura e Aplicada} \\
\mbox{} \vspace{-0.670cm} \\
{\small Universidade Federal do Rio Grande do Sul} \\
\mbox{} \vspace{-0.670cm} \\
{\small Porto Alegre, RS 91509-900, Brazil} \\
%
%
%
\nl
\mbox{} \vspace{-0.400cm} \\
{\bf Abstract} \\
\mbox{} \vspace{-0.525cm} \\
\begin{minipage}[t]{11.750cm}
{\small
\mbox{} \hspace{+0.150cm}
We provide a detailed derivation
of several fundamental properties
of bounded weak solutions
to initial value problems
for general conservative
2nd-order parabolic equations
with $p$-\mbox{\footnotesize L}aplacian
diffusion
and arbitrary
initial data
$ \;\!u_{0} \!\!\;\!\;\!\in
L^{1}(\mathbb{R}^{n}) \cap
\!\;\!\;\!L^{\infty}(\mathbb{R}^{n}) $. \\
}
\end{minipage}
\end{center}
%
%
%
%
\mbox{} \vspace{-0.650cm} \\
\setcounter{page}{1}
\mbox{} \vspace{-0.250cm} \\
%
%
%
%

%
%

%
{\bf 1. Introduction} \\
\mbox{} \vspace{-0.650cm} \\

In this work,
we provide a detailed derivation
of several fundamental properties
of (bounded, weak) solutions
of the initial value problem
for evolution
$p$-\mbox{\small L}aplacian
equations
of the type \\
\mbox{} \vspace{-0.650cm} \\
\begin{equation}
\tag{1.1$a$}
u_t \,+\: \mbox{div}\;\!\mbox{\boldmath $f$}(x,t,u)
\,+\: \mbox{div}\,\mbox{\boldmath $g$}(t,u)
\:=\;
\mu(t) \;\!\;\! \mbox{div}\,(\;\! |\;\!\nabla u\,|^{\:\!p-2}
\,\nabla u \;\!),
\end{equation}
\mbox{} \vspace{-0.850cm} \\
\begin{equation}
\tag{1.1$b$}
u(\cdot,0) \,=\,
u_{0} \in L^{1}(\mathbb{R}^{n}) \cap
L^{\infty}(\mathbb{R}^{n}).
\end{equation}
\mbox{} \vspace{-0.250cm} \\
Here,
$ p > 2 $ is constant,
$ \mu \in C^{0}(\:\![\;\!0, \infty)\,\!) $
is assumed to be positive everywhere,
and
$ \mbox{\boldmath $f$}\!\:\!= (\:\!f_{\mbox{}_{1}} \!\:\!,
f_{\mbox{}_{2}} \!\;\!,\!\;\!...\:\!, f_{\mbox{}_{\scriptstyle n}}) $,
$ \mbox{\boldmath $g$}\!\:\!= (\;\!g_{\mbox{}_{1}} \!\!\;\!\;\!,
g_{\mbox{}_{2}} \!\;\!,\!\;\!...\:\!, g_{\mbox{}_{\scriptstyle n}}) $
%
%
are
given continuous fields
such that
$ \mbox{\boldmath $g$}(t,0) = {\bf 0} $
for all $\:\!t \geq 0 \:\!$
and with
$ \!\;\!\mbox{\boldmath $f$} \!\;\!$
satisfying
the growth condition \\
\mbox{} \vspace{-0.675cm} \\
\begin{equation}
\tag{1.2}
|\,\mbox{\boldmath $f$}(x,t,\mbox{u}) \,|
\,\leq\,
\mbox{\small $F$}(t) \,|\,\mbox{u}\,|^{\:\!\kappa\,+\,1}
\quad \;\;\,
\forall \; x \in \mathbb{R}^{n} \!,
\; t \geq 0,
\; \mbox{u} \in \mathbb{R}
\end{equation}
\mbox{} \vspace{-0.250cm} \\
for some
$ \mbox{\small $F$} \in C^{0}(\!\;\!\;\![\;\!0, \infty)\,\!) $
and some constant $ \:\!\kappa \geq 0 $,
where
$ \:\!|\!\;\!\;\!\cdot\!\;\!\;\!| $
denotes the absolute value (in case of scalars)
or the Euclidean norm (in case of vectors),
as in~(1.1$a$). \linebreak
%
%
%
By a (bounded) {\em solution\/}
of (1.1) in some time interval
$ [\;\!0, \;\!\mbox{\small $T$}_{\!\ast}) $
we mean any
function \linebreak
$ u(\cdot,t) \in
C^{0}([\;\!0, \;\!\mbox{\small $T$}_{\!\ast}),
\!\;\!\;\!
L^{1}_{\mbox{\scriptsize loc}}(\mathbb{R}^{n}))
\!\;\!\;\!\cap\!\;\!\;\!
L^{p}_{\mbox{\scriptsize loc}}((\!\;\!\;\!0, \;\!\mbox{\small $T$}_{\!\ast}),
\!\;\!\;\!
W^{1,\,p}_{\mbox{\scriptsize loc}}(\mathbb{R}^{n}))
$
satisfying
the equation (1.1$a$) \linebreak
in
$ {\cal D}^{\;\!\prime}(\;\!\mathbb{R}^{n} \!\times\!\!\;\!\;\!
(\;\!0, \;\!\mbox{\small $T$}_{\!\ast})\,\!) $
with
$ {\displaystyle
\:\!
u(\cdot,0) =\:\! u_0
\,\!
} $
and
$ {\displaystyle
\;\!
u(\cdot,t) \in
L^{\infty}_{\mbox{\scriptsize loc}}([\;\!0, \;\!\mbox{\small $T$}_{\!\ast}),
\!\;\!\;\!
L^{1}(\mathbb{R}^{n}) \cap
L^{\infty}(\mathbb{R}^{n}))
} $
--- that is,
for every
$ \;\!0 < \mbox{\small $T$} \!\:\!< \mbox{\small $T$}_{\!\ast} $
given,
we have \\
\mbox{} \vspace{-0.575cm} \\
\begin{equation}
\tag{1.3$a$}
\|\, u(\cdot,t) \,
\|_{\mbox{}_{\scriptstyle L^{1}(\mathbb{R}^{n})}}
\,\leq\;
\mbox{\small $M$}_{\mbox{}_{\!1}}\!\;\!(\:\!\mbox{\small $T$}),
\quad \;\;\,
\forall \;\,
0 \:\!\leq\:\! t \:\!\leq\:\! \mbox{\small $T$}\!\;\!,
\end{equation}
\mbox{} \vspace{-0.850cm} \\
\begin{equation}
\tag{1.3$b$}
\|\, u(\cdot,t) \,
\|_{\mbox{}_{\scriptstyle L^{\infty}(\mathbb{R}^{n})}}
\leq\,
\mbox{\small $M$}_{\mbox{}_{\!\infty}}\!\;\!(\:\!\mbox{\small $T$}),
\quad \;\;
\forall \;\,
0 \:\!\leq\:\! t \:\!\leq\:\! \mbox{\small $T$} \!\;\!,
\end{equation}
\mbox{} \vspace{-0.165cm} \\
for suitable bounds
$ {\displaystyle
\mbox{\small $M$}_{\mbox{}_{\!1}}\!\;\!(\:\!\mbox{\small $T$}),
\;\!
\mbox{\small $M$}_{\mbox{}_{\!\infty}}\!\;\!(\:\!\mbox{\small $T$})
} $
depending on $ \mbox{\small $T$} $
\!\!\;\!\;\!(and the solution $u$ considered).
For the {\em local\/} (in time)
existence
of such solutions,
see e.g.\;\cite{Kalashnikov1987, Lions1969, 
WuZhaoYinLi2001, Zhao1993, Zhou2000},
while,
for global existence,
\cite{ChagasGuidolinZingano2017, Kalashnikov1987}
can be consulted.
Our main objective in this work
is to provide a complete, rigorous derivation
of important fundamental properties
possessed by the solutions,
following the lines of
\cite{BrazMeloZingano2015, BrazSchutzZingano2013, %
ChagasGuidolinZingano2017, DiBenedetto1993, %
Guidolin2015, Kalashnikov1987, WuZhaoYinLi2001}.
Thus,
for example,
in Section 2
we show that \\
\mbox{} \vspace{-0.600cm} \\
\begin{equation}
\tag{1.4}
\int_{\mbox{}_{\scriptstyle 0}}^{\;\!T}
\!\!\!\;\!
\int_{\mbox{}_{\scriptstyle \!\;\!\mathbb{R}^{n}}}
\!\!\:\!
|\, \nabla u(x,t) \,|^{\:\!p}
\:
dx\, dt
\:<\,
\infty
\end{equation}
\mbox{} \vspace{-0.100cm} \\
for every
$ {\displaystyle
\:\!
0 < \mbox{\small $T$} <
\mbox{\small $T$}_{\!\ast} \!\;\!
} $,
so that
$ {\displaystyle
\:\!
u(\cdot,t)
\in
L^{p}_{\mbox{\scriptsize loc}}
(\:\![\;\!0, \;\!\mbox{\small $T$}_{\!\ast}),
\!\;\!\;\!
W^{1,\,p}(\mathbb{R}^{n})\,\!)
} $,
along with the
monotonicity of
$ {\displaystyle
\;\!
\|\, u(\cdot,t) \,
\|_{\mbox{}_{\scriptstyle L^{1}(\mathbb{R}^{n})}}
} $
and other basic results.
In Section 3,
solutions are shown
to contract in $L^{1}(\mathbb{R}^{n})$,
so that we have \\
\mbox{} \vspace{-0.600cm} \\
\begin{equation}
\tag{1.5}
\|\, u(\cdot,t) - v(\cdot,t) \,
\|_{\mbox{}_{\scriptstyle L^{1}(\mathbb{R}^{n})}}
\;\!\leq\;
\|\, u(\cdot,0) - v(\cdot,0) \,
\|_{\mbox{}_{\scriptstyle L^{1}(\mathbb{R}^{n})}}
\end{equation}
\mbox{} \vspace{-0.225cm} \\
for any given solution pair $ \;\!u, \;\!v $,
and any $ t > 0 $
for which both solutions are defined,
provided that
the flux functions
$ \mbox{\boldmath $f$}\!, \;\!\mbox{\boldmath $g$} $
in the equation (1.1$a$) above
satisfy additional conditions,
which include \\
\mbox{} \vspace{-0.600cm} \\
\begin{equation}
\tag{1.6}
|\, \mbox{\boldmath $f$}(x,t,\mbox{u}) -
\mbox{\boldmath $f$}(x,t,\mbox{v}) \,|
\;\leq\:
\mbox{\small $K$}_{\!f}(\:\!\mbox{\small $M$}\!\:\!,
\;\!\mbox{\small $T$}) \,
|\:\mbox{u} - \mbox{v} \:|^{\:\!1 \,-\, \frac{\scriptstyle 1}{\scriptstyle p}}
\end{equation}
\mbox{} \vspace{-0.775cm} \\
\begin{equation}
\tag{1.7}
|\, \mbox{\boldmath $g$}(t,\mbox{u}) -
\mbox{\boldmath $g$}(t,\mbox{v}) \,|
\;\leq\:
\mbox{\small $K$}_{\!\:\!g}(\:\!\mbox{\small $M$}\!\:\!,
\;\!\mbox{\small $T$}) \,
|\:\mbox{u} - \mbox{v} \:|^{\:\!1 \,-\, \frac{\scriptstyle 1}{\scriptstyle p}}
\end{equation}
\mbox{} \vspace{-0.150cm} \\
for all
$ \;\!x \in \mathbb{R}^{n} \!\;\!$,
$ \;\!0 \leq t \leq \mbox{\small $T$} \!\;\!$,
$ \;\!|\,\mbox{u}\,| \leq \mbox{\small $M$}\!\:\!$,
$ \;\!|\,\mbox{v}\,| \leq \mbox{\small $M$}\!\:\!$,
for each given
$ \mbox{\small $M$} > 0 $,
$ \mbox{\small $T$} > 0 $,
where the
Lipschitz constants
$ {\displaystyle
\mbox{\small $K$}_{\!f}(\:\!\mbox{\small $M$}\!\:\!,
\;\!\mbox{\small $T$}),
\;\!
\mbox{\small $K$}_{\!\:\!g}(\:\!\mbox{\small $M$}\!\:\!,
\;\!\mbox{\small $T$})
} $
may depend
upon the values
of
$ \:\!\mbox{\small $M$} \!\;\!$,
$ \mbox{\small $T$} $ \linebreak
(see Section 3 for further details).
Also,
under such extra assumptions,
the solutions are shown to obey
a familiar
comparison principle,
as expected
for 2nd-order parabolic problems.
From this,
it follows
in particular
that
solutions are uniquely defined
by their initial data,
which is not necessarily the situation
in Section 2.

\newpage
\mbox{} \vspace{-1.200cm} \\
%

%
%
%
{\bf 2. Some fundamental basic properties} \\
\mbox{} \vspace{-0.600cm} \\
We begin by recalling
an important regularization technique
\cite{DiBenedetto1993, Urbano2008, WuZhaoYinLi2001}:
given an interval
$ I \subseteq \mathbb{R} $
(arbitrary),
$ \:\!h > 0 $ (small),
and some function
$ \:\!v(\cdot,t) \in L^{r}(\:\!I\!\;\!, \!\;\!\;\!
L^{q}_{\mbox{\scriptsize loc}}(\!\;\!\;\!\mathbb{R}^{n})) $,
where
$ \:\!q, \:\!r \in [\;\!1, \infty\:\!]$,
let
$ v_{h}(\cdot,t)
\in C^{0}(\:\!I\!\:\!, \:\!
L^{q}_{\mbox{\scriptsize loc}}(\:\!\mathbb{R}^{n})) $
be the {\em Steklov average\/} \\
\mbox{} \vspace{-0.375cm} \\
\begin{equation}
\tag{2.1}
v_{\mbox{}_{\scriptstyle h}}\!\:\!(\cdot,t)
\,:=\;
\frac{1}{h} \!\:\!
\int_{\,\!t}^{\;\!t \:\!+\:\! h}
\hspace{-0.500cm}
\tilde{v}(\cdot,\tau) \,d\tau,
\qquad
t \in I \!\:\!,
\end{equation}
\mbox{} \vspace{-0.075cm} \\
where
$ \!\;\!\;\!\tilde{v}(\cdot,\tau) = v(\cdot,\tau) $
if $ \tau \!\in I \!\;\!$,
$ \!\;\!\;\!\tilde{v}(\cdot,\tau) = 0 \!\;\!\;\!$
if $ \tau \!\notin I \!\;\!$.
For
$ {\displaystyle
\;\!
u(\cdot,t) \in
C^{0}([\;\!0, \;\!\mbox{\small $T$}_{\!\ast}),
\!\;\!\;\!
L^{1}_{\mbox{\scriptsize loc}}(\mathbb{R}^{n}))
} $
$ {\displaystyle
\cap \;\!
L^{p}_{\mbox{\scriptsize loc}}((\:\!0, \;\!\mbox{\small $T$}_{\!\ast}),
\!\;\!\;\!
W^{1,\,p}_{\mbox{\scriptsize loc}}(\mathbb{R}^{n}))
} $
solution of (1.1),
we then obtain
(see \cite{DiBenedetto1993}, Ch.\;II\:\!;
\cite{WuZhaoYinLi2001}, Ch.\:1)
that,
for any ball
$ \;\!\mbox{\small $B$}_{\mbox{}_{\!R}} \!\:\!=
\{\;\! x \in \mathbb{R}^{n} \!\!\;\!: \;\!
|\;\!x\;\!| <\!\;\! \mbox{\small $R$} \;\!\}
\:\!$: \\
\mbox{} \vspace{+0.150cm} \\
\mbox{} \hspace{+1.000cm}
$ {\displaystyle
\int_{\mbox{}_{\scriptstyle \!\:\!B_{\mbox{}_{\!\:\!R}}}}
\hspace{-0.150cm}
\Bigl\{\,
u_{\mbox{}_{\scriptstyle h,\,t}}\!\:\!(x,t) \, \phi(x)
\;\!+\;\!
\langle \;\!
\bigl[\, \mu(t) \, |\;\! \nabla u\,|^{\:\!p-2} \, \nabla u \,
\bigr]_{\!\;\!h}, \nabla \phi
\;\!\;\!\rangle
\, \Bigr\}
\;\!\;\! dx
\;=
} $ \\
\mbox{} \vspace{-0.725cm} \\
\mbox{} \hfill (2.2) \\
\mbox{} \vspace{-0.425cm} \\
\mbox{} \hspace{+3.250cm}
$ {\displaystyle
=
\int_{\mbox{}_{\scriptstyle \!\:\!B_{\mbox{}_{\!\:\!R}}}}
\hspace{-0.150cm}
\Bigl\{\;\!
\langle \,
\bigl[\;\! \mbox{\boldmath $f$}(x,t,u) \,
\bigr]_{\!\;\!h},
\:\!\nabla \phi
\;\!\;\!\rangle
\,+\;\!
\langle \,
\bigl[\;\! \mbox{\boldmath $g$}(t,u) \,
\bigr]_{\!\;\!h},
\:\!\nabla \phi
\;\!\;\!\rangle
\;\!\Bigr\}
\;\!\;\!dx
} $ \\
\mbox{} \vspace{+0.300cm} \\
for all
$ \;\!0 < t < \mbox{\small $T$}_{\!\ast} \!\;\! - h $,
\:\!and
any
$ \:\!\phi \in W^{1,\,p}_{0}(\mbox{\small $B$}_{\mbox{}_{\!\;\!R}})
\cap L^{\infty}(\mbox{\small $B$}_{\mbox{}_{\!\;\!R}}) $,
where
$ \;\!
u_{\mbox{}_{\scriptstyle h,\,t}}\!\;\!(\cdot,t)
=\!\;\!\;\!
\mbox{\footnotesize $ {\displaystyle
\frac{\partial}{\partial \:\!t} }$}
\, u_{\mbox{}_{\scriptstyle h}}
\!(\cdot,t) $
$ = \!\;\![\;\!u(\cdot,t+h) - u(\cdot,t)\:\!]
\!\;\!\;\!/\!\;\!\;\!h $
is the strong pointwise derivative
of
$ u_{\mbox{}_{\scriptstyle h}}
\!(\cdot,t) $
in
$ \mbox{\small $L$}^{\!\;\!1}(\mbox{\footnotesize $B$}_{\mbox{}_{\!R}}\!\;\!) $,
%
%
and where
$ \langle\,\cdot\;\!,\:\!\cdot\,\rangle $
denotes the
standard inner product
of a pair of $n$-dimensional vectors.
As in \cite{DiBenedetto1993, Urbano2008, WuZhaoYinLi2001},
the expression (2.2) is a very useful starting point
for the derivation
of a number of important solution properties,
as illustrated by the following results. \linebreak
\nl
\mbox{} \vspace{-0.500cm} \\
%
%
%
%
%
{\bf Proposition 2.1.}
\textit{%
Let
$ {\displaystyle
\;\!
u(\cdot,t) \in
C^{0}(\,\![\,0, \;\!\mbox{\small $T$}_{\!\ast}),
\!\;\!\;\!
L^{1}_{\mbox{\scriptsize \em loc}}(\mathbb{R}^{n})\,\!)
\!\;\!\;\!\cap\!\;\!\;\!
L^{p}_{\mbox{\scriptsize \em loc}}(\,\!
(\!\;\!\;\!0, \;\!\mbox{\small $T$}_{\!\ast}),
\!\;\!\;\!
W^{1,\,p}_{\mbox{\scriptsize \em loc}}(\mathbb{R}^{n})\,\!)
\;\cap
} $ \\
$ {\displaystyle
L^{\infty}_{\mbox{\scriptsize \em loc}}
(\:\![\;\!0, \;\!\mbox{\small $T$}_{\!\ast}),
\!\;\!\;\!
L^{1}(\mathbb{R}^{n}) \cap
L^{\infty}(\mathbb{R}^{n}))
} $
be any
given solution to the problem
$\;\!(1.1)$, $(1.2)$, \\
where
$ \,\kappa \geq 0 $.
Then
} \\
\mbox{} \vspace{-0.950cm} \\
\begin{equation}
\tag{2.3}
\int_{\mbox{}_{\scriptstyle 0}}^{\;\!T}
\!\!\!\;\!
\int_{\mbox{}_{\scriptstyle \!\;\!\mathbb{R}^{n}}}
\!\!\:\!
|\, \nabla u(x,t) \,|^{\:\!p}
\;
dx\, dt
\,<\, \infty
\end{equation}
\mbox{} \vspace{-0.075cm} \\
\textit{%
for every
$ \;\!0 < \mbox{\small $T$} \!\;\!< \mbox{\small $T$}_{\!\ast} $,
so that
$ {\displaystyle
\;\!
u(\cdot,t) \in
L^{p}_{\mbox{\scriptsize \em loc}}(\:\!
[\;\!0, \;\!\mbox{\small $T$}_{\!\ast}),
\!\;\!\;\!
W^{1,\,p}(\mathbb{R}^{n})\,\!)
} $. \\
}
%
\nl
%
%
%
%
%
\nl
%
%
\mbox{} \vspace{-0.950cm} \\
{\small
{\bf Proof.}
Let $ 0 < t_0 < \mbox{\footnotesize $T$} $.
\!Given $ \mbox{\footnotesize $R$} > 0 $,
$ \epsilon > 0 $,
let
$ \zeta_{\mbox{}_{R, \;\!{\scriptstyle \epsilon}}}
\!\in C^{2}(\mathbb{R}^{n}) $
be the cut-off function \\
\mbox{} \vspace{-0.450cm} \\
\begin{equation}
\notag
\zeta_{\mbox{}_{R, \;\!{\scriptstyle \epsilon}}}(x)
\;=\;\!\;\!
\Bigl\{\;\!
e^{\!-\, \epsilon \;\!
\,\sqrt{\:\!1 \,+\, |\;\!x\;\!|^{2}\;\!}}
-\;\!\;\!
e^{\!-\, \epsilon\;\!\sqrt{\;\!1 \;\!+\;\!R^{\:\!2}\:\!}}
\;\!\Bigr\}^{\!\:\!p}
\qquad \,
\mbox{if } \;
|\;\!x\;\!| \;\!<\:\!\mbox{\footnotesize $R$}
\end{equation}
\mbox{} \vspace{-0.125cm} \\
and
$ {\displaystyle
\;\!
\zeta_{\mbox{}_{R, \;\!{\scriptstyle \epsilon}}}(x)
\;\!=\;\!0
\;\!
} $
if
$ \;\!|\;\!x\;\!| \geq \mbox{\footnotesize $R$} $.
Taking
$ {\displaystyle
\;\!
\phi(x) \:\!=\;\!
2 \, u_{\mbox{}_{\scriptstyle h}}
\!(x,t)
\,\zeta_{\mbox{}_{R, \;\!{\scriptstyle \epsilon}}}\!\;\!(x)
} $
in (2.2) above,
integrating the resulting equation
in
$ (\!\;\!\!\;\!\;\!t_0, \!\;\!\;\!\mbox{\footnotesize $T$}) $,
and letting $ h \;\!{\scriptstyle \searrow}\;\!0 $,
%
%
we get,
letting (as always)
$ {\displaystyle
\!\;\!
\mbox{\footnotesize $B$}_{\mbox{}_{\!\:\!R}}
\!\!\;\!\;\!
} $
denote the ball
$ {\displaystyle
\!\;\!\;\!
\bigl\{\;\! x \in \mathbb{R}^{n} \!\!\;\!:
\!\;\!\;\!
|\;\!x\;\!| < \mbox{\footnotesize $R$}
\;\!\bigr\}
} $,
and setting
$ {\displaystyle
\!\;\!\;\!
\mbox{\small \boldmath $ \tilde{f} $}
\!\!\;\!\;\!:=
\mbox{\small \boldmath $ f $} + \mbox{\small \boldmath $ g $}
} $\:\!: \\
\newpage
\mbox{} \vspace{-0.950cm} \\
\mbox{} \hspace{+2.250cm}
$ {\displaystyle
\int_{\mbox{}_{\scriptstyle \!\:\!B_{\mbox{}_{\!\:\!R}}}}
\!\!\!\!\;\!
u(x,\mbox{\footnotesize$T$})^{2}
\,
\zeta_{\mbox{}_{R, \;\!{\scriptstyle \epsilon}}}\!\;\!(x)
\;\!\;\!dx
\;+\;
2 \!
\int_{\mbox{}_{\scriptstyle \!\:\!t_{\mbox{}_{0}}}}^{\:\!T}
\!\!\!
\mu(t)
\!
\int_{\mbox{}_{\scriptstyle \!\:\!B_{\mbox{}_{\!\:\!R}}}}
\!\!\!\!\;\!
|\,\nabla u\,|^{\:\!p}
\,
\zeta_{\mbox{}_{R, \;\!{\scriptstyle \epsilon}}}\!\;\!(x)
\;\!\;\!
dx\,dt
\;\!\;\!\;\!=
} $ \\
\mbox{} \vspace{+0.200cm} \\
%
\mbox{} \hspace{+0.950cm}
$ {\displaystyle
=\;\!
\int_{\mbox{}_{\scriptstyle \!\:\!B_{\mbox{}_{\!\:\!R}}}}
\!\!\!\!\;\!
u(x,t_{0})^{2}
\,
\zeta_{\mbox{}_{R, \;\!{\scriptstyle \epsilon}}}\!\;\!(x)
\;\!\;\!dx
\;-\;
2 \!\:\!
\int_{\mbox{}_{\scriptstyle \!\:\!t_{\mbox{}_{0}}}}^{\:\!T}
\!\!\!
\mu(t)
\!
\int_{\mbox{}_{\scriptstyle \!\:\!B_{\mbox{}_{\!\:\!R}}}}
\!\!\!\!\;\!
u(x,t) \,
|\,\nabla u\,|^{\:\!p - 2}
\;\!
\langle \,
\nabla u, \!\;\!\;\!
\nabla \zeta_{\mbox{}_{R, \;\!{\scriptstyle \epsilon}}}\!\;\!(x)
\;\!\rangle
\;\!\;\!
dx\,dt
} $ \\
\mbox{} \vspace{+0.050cm} \\
%
\mbox{} \hfill
$ {\displaystyle
+\;\!\;\!\;\!
2 \!\:\!
\int_{\mbox{}_{\scriptstyle \!\:\!t_{\mbox{}_{0}}}}^{\:\!T}
\!\!
\int_{\mbox{}_{\scriptstyle \!\:\!B_{\mbox{}_{\!\:\!R}}}}
\!\!\!\!\;\!
\langle \,
\mbox{\boldmath $\tilde{f}$}(x,t,u), \!\;\!\;\!
\nabla u
\;\!\rangle
\,
\zeta_{\mbox{}_{R, \;\!{\scriptstyle \epsilon}}}\!\;\!(x)
\;\!\;\!
dx\,dt
\;+\;
2 \!\:\!
\int_{\mbox{}_{\scriptstyle \!\:\!t_{\mbox{}_{0}}}}^{\:\!T}
\!\!\!
\int_{\mbox{}_{\scriptstyle \!\:\!B_{\mbox{}_{\!\:\!R}}}}
\!\!\!\!\;\!
u(x,t)
\;\!
\langle \,
\mbox{\boldmath $\tilde{f}$}(x,t,u), \!\;\!\;\!
\nabla \zeta_{\mbox{}_{R, \;\!{\scriptstyle \epsilon}}}\!\;\!(x)
\;\!\rangle
\;\!\;\!
dx\,dt
} $ \\
\mbox{} \vspace{+0.200cm} \\
\mbox{} \hfill
$ {\displaystyle
\leq\;\!\;\!
\mbox{\footnotesize $M$}_{\!\:\!1}
(\mbox{\footnotesize $T$})
\,
\mbox{\footnotesize $M$}_{\!\;\!\infty}
(\mbox{\footnotesize $T$})
\!\;\!\;\!+\!\!\;\!\;\!
\int_{\mbox{}_{\scriptstyle \!\:\!t_{\mbox{}_{0}}}}^{\:\!T}
\!\!\!
\mu(t)
\!
\int_{\mbox{}_{\scriptstyle \!\:\!B_{\mbox{}_{\!\:\!R}}}}
\!\!\!\!\;\!
|\,\nabla u\,|^{\:\!p}
\;\!
\zeta_{\mbox{}_{R, \;\!{\scriptstyle \epsilon}}}\!\;\!(x)
\;\!\;\!
dx\,dt
\;\!\;\!+\;\!\;\!
\mbox{\small %
$ {\displaystyle \frac{\,2^{\:\!p}}{p} }$}
\!
\int_{\mbox{}_{\scriptstyle \!\:\!t_{\mbox{}_{0}}}}^{\:\!T}
\!\!\!
\mu(t)
\!\!\;\!
\int_{\mbox{}_{\scriptstyle \!\:\!B_{\mbox{}_{\!\:\!R}}}}
\!\!\!\!\;\!
|\,u\,|^{\:\!p}
\,
\frac{\,
|\;\!\nabla \zeta_{\mbox{}_{R, \;\!{\scriptstyle \epsilon}}}\;\!|^{\:\!p}}
{\;\mbox{\small $ {\displaystyle %
   \zeta_{\mbox{}_{\!\:\!R, \;\!{\scriptstyle \epsilon}}}
        ^{\;\!p \;\!-\:\! 1} }$}}
\;\!\;\!
dx\,dt
} $ \\
\mbox{} \vspace{+0.050cm} \\
\mbox{} \hfill
$ {\displaystyle
+\;\;\!
2 \!\:\!
\int_{\mbox{}_{\scriptstyle \!\:\!t_{\mbox{}_{0}}}}^{\:\!T}
\!\!\!
\mbox{\footnotesize $F$}(t)^{\mbox{}^{\scriptstyle \!\;\!
\frac{\scriptstyle p}{\scriptstyle p \;\!-\:\! 1} }}
\!\;\!
\mu(t)^{\mbox{}^{\scriptstyle \!\!\!
-\, \frac{\scriptstyle 1}{\scriptstyle p\;\!-\:\!1} }}
\!\!\!\!\;\!\;\!
\int_{\mbox{}_{\scriptstyle \!\:\!B_{\mbox{}_{\!\:\!R}}}}
\!\!\!\!\;\!
|\,u\,|^{\mbox{}^{\scriptstyle \!\;\!
(1 \:\!+\;\! \kappa) \,
\frac{\scriptstyle p}{\scriptstyle p \;\!-\:\!1} }}
\!\;\!\;\!
\zeta_{\mbox{}_{R, \;\!{\scriptstyle \epsilon}}}
\!\;\!\;\!
dx\,dt
%
%
\;\!\;\!+\;\!\!\;\!\;\!
2 \!\!\;\!
\int_{\mbox{}_{\scriptstyle \!\:\!t_{\mbox{}_{0}}}}^{\:\!T}
\!\!\!\!\;\!
\mbox{\footnotesize $F$}(t)
\!
\int_{\mbox{}_{\scriptstyle \!\:\!B_{\mbox{}_{\!\:\!R}}}}
\!\!\!\!\;\!
|\,u\,|^{\mbox{}^{\scriptstyle \:\!2 \;\!+\;\! \kappa }}
\,
|\, \nabla \zeta_{\mbox{}_{R, \;\!{\scriptstyle \epsilon}}} \;\!|
\;\!\;\!
dx\,dt
} $ \\
\mbox{} \vspace{+0.050cm} \\
%
\mbox{} \hfill
$ {\displaystyle
+\;\,
4 \;
\mbox{\footnotesize $G$}(\mbox{\footnotesize $T$})
\!
\int_{\mbox{}_{\scriptstyle \!\:\!t_{\mbox{}_{0}}}}^{\:\!T}
\!\!
\int_{\mbox{}_{\scriptstyle \!\:\!B_{\mbox{}_{\!\:\!R}}}}
\!\!\!\!\;\!\;\!
|\, u(x,t)\,|
\;
|\, \nabla \zeta_{\mbox{}_{R, \;\!{\scriptstyle \epsilon}}} \;\!|
\;
dx\,dt
} $ \\
\mbox{} \vspace{+0.250cm} \\
by (1.2), (1.3)
and
Young's inequality
(see e.g.\;\cite{Evans1998}, p.\;622),
where
$ \;\!\mbox{\footnotesize $M$}_{\mbox{}_{\!1}} \!\;\!$,
$ \!\;\!\mbox{\footnotesize $M$}_{\mbox{}_{\!\infty}}\!\;\!$
are given in (1.3)
and
$ {\displaystyle
\!\;\!\;\!
\mbox{\footnotesize $G$}(\mbox{\footnotesize $T$})
\!\;\!\;\!=\;\!
\sup\;
\{\;\!\;\!
|\;\!\;\!\mbox{\boldmath $g$}(t,\mbox{\small v})\,|
\!\!\;\!\;\! : \,
0 < t < \mbox{\footnotesize $T$}, \;
|\,\mbox{\small v}\,| \!\;\!\;\!<\!\;\!\;\!
\mbox{\footnotesize $M$}_{\mbox{}_{\!\infty}}\!\;\!
(\mbox{\footnotesize $T$}) \;\!\}
} $.
Letting
$ \mbox{\footnotesize $R$} \;\!\mbox{\scriptsize $\nearrow$}\, \infty $,
$ \epsilon \,\mbox{\scriptsize $\searrow$}\, 0 $
and
$ t_0 \!\;\!\;\!\mbox{\scriptsize $\searrow$}\, 0 $ \linebreak
(in this order),
we then
obtain,
by (1.3$a$)
and since
$ {\displaystyle
\;\!
|\, \nabla \zeta_{\mbox{}_{\!\;\!R, \;\!{\scriptstyle \epsilon}}} \:\!|^{\:\!p}
/\;\!\zeta_{\mbox{}_{\!\;\!R, \;\!{\scriptstyle \epsilon}}}^{\;\!p \;\!-\:\!1}
\!\;\!\leq
\bigl(\:\!p \,\epsilon
\:\!\bigr)^{\scriptstyle p} \;\!
e^{\scriptstyle -\, p\;\! \epsilon\;\!
\sqrt{\:\! 1 \,+\, |\;\!x\;\!|^{\:\!2}\;\!}}
\!
} $, \\
\mbox{} \vspace{-0.001cm} \\
\mbox{} \hfill
$ {\displaystyle
\|\, u(\cdot,\mbox{\footnotesize $T$}) \,
\|_{\mbox{}_{\scriptstyle L^{2}(\mathbb{R}^{n})}}
  ^{\:\!2}
+\!\;\!
\int_{\mbox{}_{\scriptstyle \!\:\!0}}
    ^{\:\!T}
\!\!\!
\mu(t)
\!
\int_{\mbox{}_{\scriptstyle \!\;\!\mathbb{R}^{n}}}
\!\!\!
|\,\nabla u \,|^{\:\!p}
\,dx\,dt
\;\leq\;
\mbox{\footnotesize $M$}_{\mbox{}_{\!\infty}}
\!\!\;\!\;\!(\mbox{\footnotesize $T$})
\,
\|\, u_{0} \;\!
\|_{\mbox{}_{\scriptstyle L^{1}(\mathbb{R}^{n})}}
\!\;\!\:\!+\!\!\;\!\;\!
\int_{\mbox{}_{\scriptstyle \!\:\!0}}^{\:\!T}
\!\!\!\!\;\!
\mbox{w}(t)
\,
\|\, u(\cdot,t) \,
\|_{{\scriptstyle L^{q^{\prime}}\!(\mathbb{R}^{n})}}
  ^{\:\!q^{\prime}}
dt
} $ \\
\mbox{} \vspace{+0.075cm} \\
where
$ {\displaystyle
\,
\mbox{w}(t)
\!\;\!\;\!=\!\;\!\;\!
2 \;\!\!\;\!\;\!
\mbox{\footnotesize $F$}(t)^{\mbox{}^{\scriptstyle \!\;\!
\frac{\scriptstyle p}{\scriptstyle p \;\!-\:\! 1} }}
\!\!\;\!\;\!
\mu(t)^{\mbox{}^{\scriptstyle \!\!\!
-\, \frac{\scriptstyle 1}{\scriptstyle p\;\!-\:\!1} }}
\!
} $
and
$ {\displaystyle
\:\!
q^{\prime} \!\:\!=\;\!
(1 + \kappa) \;\! p \:\!/ (p - 1 )
} $.
This shows (2.3).
\mbox{} \hfill $\Box$ \\
}
%
%
\mbox{} \vspace{-0.500cm} \\

The next result gives one form of
the basic energy inequalities that
can be obtained
for weak solutions
$ {\displaystyle
\;\!
u(\cdot,t) \in
C^{0}([\;\!0, \;\!\mbox{\small $T$}_{\!\ast}),
\!\;\!\;\!
L^{1}_{\mbox{\scriptsize loc}}(\mathbb{R}^{n}))
\:\!\cap\:\!
L^{p}_{\mbox{\scriptsize loc}}
((\!\;\!\;\!0, \;\!\mbox{\small $T$}_{\!\ast}),
\:\!
W^{1,\,p}_{\mbox{\scriptsize loc}}(\mathbb{R}^{n}))
\;\!\;\!\cap
} $
$ {\displaystyle
L^{^\infty}_{\mbox{\scriptsize loc}}
([\;\!0, \;\!\mbox{\small $T$}_{\!\ast}),
\!\;\!\;\!
L^{1}(\mathbb{R}^{n})
\cap
L^{\infty}(\mathbb{R}^{n}))
} $
of problem (1.1), (1.2),
which plays a key role in
\cite{ChagasGuidolinZingano2017}.

\nl
%
%
%
%
%
{\bf Proposition 2.2.}
\textit{%
\!\!\;\!Under the same assumptions
of \,{\sc\small Proposition 2.1}\!\:\!
above,
we have,
for each $ \!\;\!\;\!q \geq 2 $,
%
%
%
%
%
%
that
$ {\displaystyle
\,
\|\, u(\cdot,t) \,
\|_{\scriptstyle L^{q}(\mathbb{R}^{n})}
  ^{\:\!q}
\!\:\!
} $
is absolutely continuous
in
$ \;\! t \in (\!\;\!\;\!0, \!\;\!\;\!\mbox{\small $T$}_{\!\ast}) $.
\!\;\!Moreover,
there exists
$ \!\;\!\;\!E_{q} \subset (\:\!0, \mbox{\small $T$}_{\!\ast}) $
with zero Lebesgue measure
such that
} \\
\mbox{} \vspace{-0.575cm} \\
\begin{equation}
\tag{2.4}
\mbox{\small $ {\displaystyle
\frac{d}{d\:\!t} }$}\;
\|\, u(\cdot,t) \,
\|_{\mbox{}_{\scriptstyle L^{2}(\mathbb{R}^{n})}}
  ^{\:\!2}
\!\;\!+\:
2 \, \mu(t)
\!\!\;\!
\int_{\mathbb{R}^{n}}
\!\!
|\, \nabla u \,|^{\:\!p}
\: dx \;
\leq\;
2 \,
\mbox{\small $F$}(t)
\!\!\;\!
\int_{\mathbb{R}^{n}}
\!\!
|\, u(x,t)\,|^{\:\!\kappa \:\!+\:\!1}
\,
|\, \nabla u \,|
\: dx
\end{equation}
\mbox{} \vspace{-0.100cm} \\
\textit{%
for all
$ {\displaystyle
\,
t \in (\:\!0, \:\!\mbox{\small $T$}_{\!\ast})
\!\;\!\setminus\!\;\! E_{2}
\!\;\!
} $
$($if $\:\! q = 2 \:\!)$,
and
} \\
\mbox{} \vspace{-0.610cm} \\
\begin{equation}
\notag
\mbox{\small $ {\displaystyle
\frac{d}{d\:\!t} }$}\;
\|\, u(\cdot,t) \,
\|_{\mbox{}_{\scriptstyle L^{q}(\mathbb{R}^{n})}}
  ^{\:\!q}
+\:
q\;\!(q - 1) \, \mu(t)
\!\!\;\!
\int_{\mathbb{R}^{n}}
\!\!
|\, u(x,t)\,|^{\:\!q \:\!-\:\! 2}
\,
|\, \nabla u \,|^{\:\!p}
\: dx
\end{equation}
\mbox{} \vspace{-0.850cm} \\
\mbox{} \hfill (2.5) \\
\mbox{} \vspace{-1.100cm} \\
\begin{equation}
\notag
\leq\;
q\;\!(q - 1) \,
\mbox{\small $F$}(t)
\!\!\;\!
\int_{\mathbb{R}^{n}}
\!\!
|\, u(x,t)\,|^{\:\!q \:\!-\:\! 1 \:\!+\:\! \kappa}
\,
|\, \nabla u \,|
\; dx
\end{equation}
\mbox{} \vspace{-0.125cm} \\
\textit{%
for all
$ {\displaystyle
\,
t \in (\:\!0, \:\!\mbox{\small $T$}_{\!\ast})
\!\;\!\setminus\!\;\! E_{q}
\!\;\!
} $
$($if $\:\!q > 2 \:\!)$,
where
$ \mbox{\small $F$}(t) \!\;\!$
is given in $(1.2)\!\;\!$
above. \\
}
%
%
%
%
{\small
{\bf Proof.}
Given $ \:\!0 < t_0 < t < \mbox{\footnotesize $T$}_{\!\!\;\!\;\!\ast} $,
$ \mbox{\footnotesize $R$} > 0 $,
let
$ \:\!\zeta_{\mbox{}_{R}}\!\;\!(x) = \zeta(x/\mbox{\footnotesize $R$}) $,
where
$ \:\!\zeta \in C^{1}(\mathbb{R}^{n}) $
is such that
$ \zeta(x) = 1 \:\!$
if $ \:\!|\:\! x \:\!| \leq 1/2 $,
%
%
$ \zeta(x) = 0 \:\!$
if $ \:\!|\:\! x \:\!| > 1 $,
$ \:\!0 \leq \zeta \leq 1 $
for all $\:\! x \in \mathbb{R}^{n} \!$.
\!\;\!We begin with
$\!\;\!\;\! q > 2 \:\!$:
Taking
$ S \in C^{2}(\mathbb{R}) $
such that
$ S^{\prime}(\mbox{u}) = -\;\!1 \:\!$
if \:\!$ \mbox{u} \leq -\;\!1 $,
$ \:\!S^{\prime}(\mbox{u}) = 1 \:\!$
if $ \:\!\mbox{u} \geq 1 $,
$ \:\!S(0) = 0 \:\!$
and
$ \:\!S^{\prime}(\mbox{u}) \geq 0 \:\!$
for all $ \:\!\mbox{u} \in \mathbb{R}$,
let
$ L(\mbox{u}) := \int_{0}^{\:\!\mbox{\scriptsize u}}
\!\;\!S(\mbox{v}) \,d\mbox{v} $,
and,
for each $ \delta > 0 $,
$ L_{\delta}(\mbox{u}) := \delta \, L(\mbox{u}/\delta) $.
(\:\!This gives
$ L_{\delta}(\mbox{u}) \rightarrow |\;\!\mbox{u}\;\!| $
as $ \delta \rightarrow 0 $,
uniformly in $ \mbox{u} \in \mathbb{R} $.)
Setting
$ {\displaystyle
\;\!
\Phi_{\delta}(\mbox{u})
:=
L_{\delta}(\mbox{u})^{q}
\!\;\!
} $,
let
us take
in (2.2)
$ {\displaystyle
\phi(x) =
\Phi_{\delta}^{\prime}(u_{\mbox{}_{\scriptstyle h}}
\!(x,t))
\,\zeta_{\mbox{}_{R}}\!\;\!(x)
} $.
\!Integrating (2.2)
in
$ (\!\;\!\!\;\!\;\!t_0, \!\;\!\!\;\!\;\!t) $
and letting
$ h \rightarrow 0 $,
$ \delta \rightarrow 0 $
and then
$ \:\!\mbox{\footnotesize $R$} \rightarrow \infty $,
we get,
by (1.3) and (2.3)
above, \\
\mbox{} \vspace{-0.600cm} \\
\begin{equation}
\notag
\|\, u(\cdot,t) \,
\|_{\mbox{}_{\scriptstyle L^{q}(\mathbb{R}^{n})}}
  ^{\:\!q}
+\:
q\,(q - 1)
\!
\int_{\mbox{}_{\scriptstyle \!\;\!t_{\mbox{}_{0}}}}^{\;\!t}
\!\!\;\!
\mu(\tau)
\!\:\!
\int_{\mbox{}_{\scriptstyle \mathbb{R}^{n}}}
\!\!
|\, u \,|^{\:\!q \;\!-\;\! 2}
\:
|\,\nabla u\,|^{\:\!p}
\:dx\,d\tau
\end{equation}
\begin{equation}
\notag
= \;
\|\, u(\cdot,t_0) \,
\|_{\mbox{}_{\scriptstyle L^{q}(\mathbb{R}^{n})}}
  ^{\:\!q}
\!\;\!+\:
q \,(q - 1)
\!
\int_{\mbox{}_{\scriptstyle t_{\mbox{}_{0}}}}^{\:\!t}
\!\;\!
\int_{\mbox{}_{\scriptstyle \mathbb{R}^{n}}}
\!\!
|\,u\,|^{\:\!q - 2} \,
\langle \,
\mbox{\boldmath $f$}(x,\tau,u), \;\! \nabla u
\,\rangle
\;dx\,d\tau,
\end{equation}
\nl
\mbox{} \vspace{-0.600cm} \\
from which the result is obtained
from (1.2), (2.3) and Lebesgue's differentiation theorem.
For the case $\:\!q = 2 \:\!$
we proceed similarly,
using
$ {\displaystyle
\:\!
\phi(x) =
u_{\mbox{}_{\scriptstyle h}}
\!(x,t)
\,\zeta_{\mbox{}_{R}}\!\;\!(x)
\:\!
} $
in (2.2)
above.
\mbox{} \hfill $\Box$ \\
}
%
%
\mbox{} \vspace{-0.625cm} \\

Sometimes
(as in \mbox{\small \sc Propositions 2.3},
\mbox{\small \sc 2.4} below)
the following extra assumption on
$ \:\!\mbox{\boldmath $g$} $
is also needed:
given any
$ \:\!\mbox{\small $T$} > 0 $,
there exists some constant
$ C(\mbox{\small $T$}) $
such that \\
\mbox{} \vspace{-0.675cm} \\
\begin{equation}
\tag{2.6}
|\;\!\;\! \mbox{\boldmath $g$}(t,\mbox{u})\,|
\;\!\;\!\leq\,
C(\mbox{\small $T$})
\;\!\;\!
|\,\mbox{u}\,|^{\mbox{}^{\scriptstyle 1 \,-\,
                \frac{\scriptstyle 1}{\scriptstyle p}}}
\qquad
\forall \;\;
|\,\mbox{u}\,| \:\!\ll\:\! 1,
\;\,
0 < t < \mbox{\small $T$}\!\;\!.
\end{equation}
\mbox{} \vspace{-0.100cm} \\
%
%
%
%
%
{\bf Proposition 2.3.}
\textit{%
\!\!\;\!Under the same assumptions
of \,{\sc\small Proposition 2.1}\!\:\!
above,
we have \linebreak
}
\mbox{} \vspace{-0.700cm} \\
\begin{equation}
\tag{2.7}
\mbox{} \hspace{+1.000cm}
\|\, u(\cdot,t) \,
\|_{\mbox{}_{\scriptstyle L^{1}(\mathbb{R}^{n})}}
\leq\;
\|\, u_0 \;\!
\|_{\mbox{}_{\scriptstyle L^{1}(\mathbb{R}^{n})}}
\!\:\!,
\qquad \;\,
\forall \;\;
0 < t < \mbox{\small $T$}_{\!\ast}
\end{equation}
\mbox{} \vspace{-0.225cm} \\
\textit{%
provided that
$\;\!($i\/$)$
$ \kappa \:\!\geq\:\! 1 - 2/p $,
\;\!or that
$\;\!($ii\/$)$
$ p \:\!\geq\:\! n \:\!$
and $\;\!(2.6)$ holds.
} \\
%
%
\mbox{} \vspace{-0.050cm} \\
%
%
{\small
{\bf Proof.}
Let
$ L_{\delta} \in C^{3}(\mathbb{R}) $,
$ \zeta_{\mbox{}_{R}} \!\in C^{1}(\mathbb{R}^{n}) $
be constructed
as in the proof of \mbox{\footnotesize \sc Proposition 2.2},
and take (2.2) with
$ {\displaystyle
\phi(x) =
L_{\delta}^{\prime}(u_{\mbox{}_{\scriptstyle h}}
\!(x,t))
\,\zeta_{\mbox{}_{R}}\!\;\!(x)
} $.
If
$ \kappa \geq 1 - 2/p $,
we may proceed as follows:
integrating (2.2) in $ (t_0, t) $
and letting
$ \:\!h \rightarrow 0 $, $ t_0 \rightarrow 0 $
and
$ \mbox{\footnotesize $R$} \rightarrow \infty $,
we obtain \\
\mbox{} \vspace{-0.600cm} \\
\begin{equation}
\notag
\int_{\mbox{}_{\scriptstyle \mathbb{R}^{n}}}
\!\!\!\!\;\!
L_{\delta}(u(x,t))\, dx
\;\:\!\leq\;
\|\, u_0 \;\!
\|_{\mbox{}_{\scriptstyle L^{1}(\mathbb{R}^{n})}}
\:\!+
\int_{\mbox{}_{\scriptstyle \!\;\!0}}^{\mbox{\footnotesize $\:\!t$}}
\!\!\;\!
\mu(\tau)^{\mbox{}^{\scriptstyle \!\! -\,
\frac{\scriptstyle 1}{\scriptstyle p\,-\;\!1}}}
\!\!\!\:\!
\int_{\mbox{}_{\scriptstyle \mathbb{R}^{n}}}
\!\!\!\!\;\!
L_{\delta}^{\prime\prime}(u(x,\tau))
\;\!\;\!
|\, \mbox{\boldmath $f$}(x,\tau,u)\,
|^{\mbox{}^{\scriptstyle \:\!\frac{\scriptstyle p}{\scriptstyle p \,-\;\!1}}}
\;\!dx
\:d\tau,
\end{equation}
\mbox{} \vspace{-0.050cm} \\
from which (2.7) is obtained
by letting $ \;\! \delta \rightarrow 0 $
(because $ \:\!(\kappa + 1)\;\!p\;\!/(p-1) \:\!\geq\:\! 2 \;\!$
in this case).
In case ({\em ii\/}),
we let instead
$ \:\! h \rightarrow 0 $,
$ t_0 \rightarrow 0 $
and $ \delta \rightarrow 0 $,
which gives,
by (2.6), \\
\mbox{} \vspace{-0.550cm} \\
\begin{equation}
\notag
\int_{\mbox{}_{\scriptstyle \!B_{\mbox{}_{\!\;\!R}}}}
\!\!\!\!\;\!
|\, u(x,t)\,| \;
\zeta_{\mbox{}_{R}}\!\;\!(x)
\;\!\;\! dx
\;\:\!\leq\;\:\!
\|\, u_0 \;\!
\|_{\mbox{}_{\scriptstyle L^{1}(\mathbb{R}^{n})}}
\;\!+
\int_{\mbox{}_{\scriptstyle \!\;\!0}}^{\mbox{\footnotesize $\:\!t$}}
\!\!\!\;\!\;\!
\mu(\tau)
\!\;\!
\int_{\mbox{}_{\scriptstyle \!B_{\mbox{}_{\!\;\!R}}}}
\!\!\!\!\;\!
|\, \nabla u(x,\tau)\,|^{\mbox{}^{\scriptstyle \:\!p \;\!-\:\! 1}}
\;\!
|\, \nabla \zeta_{\mbox{}_{R}}\!\;\!(x) \,|
\;
dx \;\!\;\!d\tau
\end{equation}
\mbox{} \vspace{-0.100cm} \\
\mbox{} \hspace{+6.700cm}
$ {\displaystyle
+\;\,
\mbox{\footnotesize $K$}\!\;\!
(\mbox{\footnotesize $M$}\!\;\!, \:\!t)
\!\:\!
\int_{\mbox{}_{\scriptstyle \!\;\!0}}^{\mbox{\footnotesize $\:\!t$}}
\!\!\;\!
\int_{\mbox{}_{\scriptstyle \!B_{\mbox{}_{\!\;\!R}}}}
\!\!\!\!\;\!
|\,u(x,\tau)\,
|^{\mbox{}^{\scriptstyle \!
\frac{\scriptstyle p \,-1}{\scriptstyle p} }}
|\, \nabla \zeta_{\mbox{}_{R}}\!\;\!(x) \,|
\;
dx \;\!\;\!d\tau
} $ \\
\nl
\mbox{} \vspace{-0.425cm} \\
for some constant
$ {\displaystyle
\:\!
\mbox{\footnotesize $K$} \!= \mbox{\footnotesize $K$}\!\;\!
(\mbox{\footnotesize $M$}\!\;\!, \:\!t)
} $
depending upon
$ \mbox{\footnotesize $M$} \!\;\!$
(the maximum size of
$ {\displaystyle
\;\!
\|\, u(\cdot,\tau) \,
\|_{\mbox{}_{\scriptstyle L^{\infty}(\mathbb{R}^{n})}}
\!\;\!
} $,
$ 0 \leq \tau \leq t \:\!$)
and $\:\!t$.
Letting
$ \mbox{\footnotesize $R$} \rightarrow \infty $,
this gives (2.7),
since we are now assuming
$ \;\! p \geq n $.
\mbox{} \hfill $\Box$ \\
%
}
\nl
%
%
%
%
{\bf Remark 2.1.}
%
%
\!In addition to conditions ({\em i\/}) and ({\em ii\/})
of \mbox{\small \sc Proposition 2.3},
if $ \mbox{\boldmath $g$} $ satisfies (2.6)
with exponent 1
(cf.\;(2.9) below),
then all solutions to (1.1), (1.2)
constructed by parabolic regularization
satisfy (2.7)
when $ \:\!p \geq 3 $:
see \cite{Guidolin2015}, Ch.\;2,
and \mbox{\small \sc Remark 2.3}. \\
%
%
%
%
\nl
%
%
%
%
{\bf Remark 2.2.}
When (2.7) is valid,
it follows more generally
that we have,
by the same argument:
$ {\displaystyle
\|\, u(\cdot,t) \,
\|_{\mbox{}_{\scriptstyle L^{1}(\mathbb{R}^{n})}}
\leq\;\!
\|\, u(\cdot,t_0) \,
\|_{\mbox{}_{\scriptstyle L^{1}(\mathbb{R}^{n})}}
\!\;\!
} $
for all
$ \;\! 0 \leq t_0 \leq t < \mbox{\small $T$}_{\!\ast} $,
\:\!so that
$ {\displaystyle
\|\, u(\cdot,t) \,
\|_{\mbox{}_{\scriptstyle L^{1}(\mathbb{R}^{n})}}
\!\;\!
} $
is then
monotonically decreasing in $\:\!t$. \\
%
%
%
%
\nl
\mbox{} \vspace{-0.450cm} \\
%
%
%
%
%
{\bf Proposition 2.4.}
\textit{%
\!Let
$ {\displaystyle
\,
u(\cdot,t) \in
C^{0}(\,\![\,0, \;\!\mbox{\small $T$}_{\!\ast}),
\!\;\!\;\!
L^{1}_{\mbox{\scriptsize \em loc}}(\mathbb{R}^{n})\,\!)
\!\:\!\cap\!\:\!
L^{\infty}_{\mbox{\scriptsize \em loc}}(\,\!
[\,0, \;\!\mbox{\small $T$}_{\!\ast}),
\!\;\!
L^{1}(\mathbb{R}^{n})
\!\:\!\cap\!\:\!
L^{\infty}(\mathbb{R}^{n})\,\!)
} $
$ {\displaystyle
\cap\;\!
L^{p}_{\mbox{\scriptsize \em loc}}(\,\!
(\!\;\!\;\!0, \;\!\mbox{\small $T$}_{\!\ast}),
\!\;\!\;\!
W^{1,\,p}_{\mbox{\scriptsize \em loc}}(\mathbb{R}^{n})\,\!)
} $
be any
solution to
$(1.1)$, $(1.2)$.
\!If $\, p \geq n $
and $(2.6)$ holds, \linebreak
%
%
then
$ {\displaystyle
\;\!
u(\cdot,t) \in
C^{0}(\,\![\,0, \;\!\mbox{\small $T$}_{\!\ast}),
\!\;\!\;\!
L^{1}(\mathbb{R}^{n})\,\!)
} $.
\!\:\!$(\!\;\!\;\!$In particular,
$ {\displaystyle
\;\!
\|\, u(\cdot,t) - u_{0} \;\!
\|_{\mbox{}_{\scriptstyle L^{1}(\mathbb{R}^{n})}}
\!\rightarrow\:\! 0
} $
as
$ \:\! t \;\!\mbox{\footnotesize $\searrow$}\;\! 0 $.$)$
Moreover,
the solution mass is conserved,
i.e., \\
}
\mbox{} \vspace{-0.600cm} \\
\begin{equation}
\tag{2.8}
\int_{\mbox{}_{\scriptstyle \mathbb{R}^{n}}}
\!\!\!\:\!
u(x,t) \: dx
\;=\;\!
\int_{\mbox{}_{\scriptstyle \mathbb{R}^{n}}}
\!\!\!\:\!
u_0(x) \, dx,
\qquad \;\,
\forall \;\,
0 \:\!<\:\! t \:\!<\:\! \mbox{\small $T$}_{\!\ast}.
\end{equation}
%
%
%
%
\mbox{} \vspace{+0.100cm} \\
%
%
%
{\small
{\bf Proof.}
We begin
by showing
that
$ {\displaystyle
\;\!
u(\cdot,t) \in
C^{0}(\,\![\,0, \;\!\mbox{\small $T$}_{\!\ast}),
\!\;\!\;\!
L^{1}(\mathbb{R}^{n})\,\!)
} $.
The following argument
is adapted from
\cite{BrazSchutzZingano2013},
Theorem 2.1.
%
%
Since
$ \;\!u(\cdot,t)\;\!$
is already known
to be continuous in
$ L^{1}_{\mbox{\scriptsize loc}}(\mathbb{R}^{n}) $,
it is
sufficient
to show that,
given
$ \;\!0 < \mbox{\footnotesize $T$} <
\mbox{\footnotesize $T$}_{\!\ast} $
arbitrary,
we have
$ {\displaystyle
\|\, u(\cdot,t) \,
\|_{\scriptstyle \:\!L^{1}(|\,x\,|\,>\,R)}
} $
uniformly small
(say, $ O(\epsilon) $)
for all $ \;\!0 < t \leq \mbox{\footnotesize $T$} $
provided that
we choose
$ {\displaystyle
\!\;\!\;\!
\mbox{\footnotesize $R$} =
\mbox{\footnotesize $R$}(\epsilon,\mbox{\footnotesize $T$})
\gg 1
} $.
Let then
$ \epsilon > 0 $,
$ \;\!0 < \mbox{\footnotesize $T$} <
\mbox{\footnotesize $T$}_{\!\ast} $
be given,
and let
$ \zeta_{\scriptstyle R,\:\!S}
\!\:\!\in C^{1}(\mathbb{R}^{n}) $
be a cut-off function
satisfying:
$ 0 \leq \zeta_{\scriptstyle R,\:\!S} \leq 1 \:\!$
everywhere,
and
$ \:\!\zeta_{\scriptstyle R,\:\!S}(x) = 0 \;\!$
if $ \:\! |\,x\,| < \mbox{\footnotesize $R$}/2 $,
$ \:\!\zeta_{\scriptstyle R,\:\!S}(x) = 1 $
if $ \:\!\mbox{\footnotesize $R$} < |\,x\,| <
\mbox{\footnotesize $R$} + \mbox{\footnotesize $S$} $,
$ \:\!\zeta_{\mbox{}_{\scriptstyle R,\:\!S}}(x) = 0 \;\!$
if $ \:\! |\,x\,| > \mbox{\footnotesize $R$} +
2 \:\!\mbox{\footnotesize $S$} $,
with
$ \:\!|\;\!\nabla \zeta_{R,\:\!S}(x)\,| \leq
\mbox{\footnotesize $C$}/\mbox{\footnotesize $R$} \;\!$
if $ \:\! |\,x\,| < \mbox{\footnotesize $R$} $
and
$ \:\!|\;\!\nabla \zeta_{R,\:\!S}(x)\,| \leq
\mbox{\footnotesize $C$}/\mbox{\footnotesize $S$} \;\!$
if $ \:\! \mbox{\footnotesize $R$} + \mbox{\footnotesize $S$}
< |\,x\,| < \mbox{\footnotesize $R$} + 2 \:\!\mbox{\footnotesize $S$} $,
for some constant
$ \mbox{\footnotesize $C$} $
independent of
$ \mbox{\footnotesize $R$}, \mbox{\footnotesize $S$} > 0 $.
%
Given $ 0 < t_0 < t \leq \mbox{\footnotesize $T$} $,
$ h > 0 $, $ \delta > 0 $,
let $ L_{\delta} \in C^{3}(\mathbb{R}) $
be the regularized absolute value function
introduced in the proof of
{\footnotesize \sc Proposition 2.2}.
Taking
$ \phi(x) = L_{\delta}^{\prime}
(u_{\mbox{}_{\scriptstyle h}}\!(x,t))
\;\! \zeta_{\mbox{}_{\scriptstyle R,\:\!S}}\!\;\!(x) $
in (2.2),
and integrating the result
in $ \,\!(\,\!t_0, \:\!t\,\!) $,
we get,
letting
$ \:\! h \rightarrow 0 $,
$ t_0 \rightarrow 0 $
and
$ \delta \rightarrow 0 $, \\
\mbox{} \vspace{-0.475cm} \\
\begin{equation}
\notag
\int_{\mbox{}_{\scriptstyle \!
R/2 \,<\, |\,x\,| \,<\, R \;\!+\;\! 2 \;\!S}}
\hspace{-2.250cm}
|\, u(x,t) \,| \,
\zeta_{\mbox{}_{\scriptstyle R, \:\!S}}(x)
\;\!\;\! dx
\;\leq\;\!
\int_{\mbox{}_{\scriptstyle \!\;\!
|\,x\,| \,>\, R/2}}
\hspace{-1.075cm}
|\, u_0(x) \,| \;\!\;\! dx
\;+\;
\mbox{\footnotesize $I$}
(\mbox{\footnotesize $R$}, \,\!\mbox{\footnotesize $S$})
\;\!\;\!+\,
\mbox{\footnotesize $J$}_{\mbox{}_{\scriptstyle \!\:\!1}}
\!\;\!(\mbox{\footnotesize $R$})
\;\!\;\!+\,
\mbox{\footnotesize $J$}_{\mbox{}_{\scriptstyle \!\,\!2}}
\!\;\!(\mbox{\footnotesize $R$}, \,\!\mbox{\footnotesize $S$})
\;\!\;\!+\,
\mbox{\footnotesize $H$}_{\mbox{}_{\scriptstyle \!\:\!1}}
\!\,\!(\mbox{\footnotesize $R$})
\;\!\;\!+\,
\mbox{\footnotesize $H$}_{\mbox{}_{\scriptstyle \!2}}
\!\;\!(\mbox{\footnotesize $R$}, \,\!\mbox{\footnotesize $S$})
\end{equation}
\mbox{} \vspace{+0.150cm} \\
by (1.2), (1.3)
and (2.3),
where \\
\mbox{} \vspace{-0.525cm} \\
\begin{equation}
\notag
\mbox{\footnotesize $I$}
(\mbox{\footnotesize $R$}, \,\!\mbox{\footnotesize $S$})
\;=\;\!
\int_{\mbox{}_{\scriptstyle 0}}^{\:\!T}
\!\!\!
\mbox{\footnotesize $F$}(\tau)
\!\!\;\!\;\!
\int_{\mbox{}_{\scriptstyle \!\:\!
R/2 \,<\, |\,x\,| \,<\, R \;\!+\;\! 2 \;\!S}}
\hspace{-2.250cm}
|\, u(x,\tau) \,|^{\:\!\kappa \;\!+\;\! 1}
\;\!\;\!
|\,\nabla \zeta_{\mbox{}_{\scriptstyle R, \:\!S}}(x)\,|
\;
dx \,d\tau,
\end{equation}
\mbox{} \vspace{-0.425cm} \\
\begin{equation}
\notag
\mbox{\footnotesize $J$}_{\mbox{}_{\scriptstyle \!\:\!1}}
\!\:\!(\mbox{\footnotesize $R$})
\;=\;\!
\int_{\mbox{}_{\scriptstyle 0}}^{\:\!T}
\!\!\!
\mu(\tau)
\!\!\;\!\;\!
\int_{\mbox{}_{\scriptstyle \!\:\!
R/2 \,<\, |\,x\,| \,<\, R}}
\hspace{-1.575cm}
|\, \nabla u \,|^{\:\!p \;\!-\;\! 1}
\;\!\;\!
|\,\nabla \zeta_{\mbox{}_{\scriptstyle R, \:\!S}}(x)\,|
\;
dx \,d\tau,
\end{equation}

\mbox{} \vspace{-1.250cm} \\
\begin{equation}
\notag
\mbox{\footnotesize $J$}_{\mbox{}_{\scriptstyle \!\:\!2}}
\!\;\!(\mbox{\footnotesize $R$},\,\!\mbox{\footnotesize $S$})
\;=\;\!
\int_{\mbox{}_{\scriptstyle 0}}^{\:\!T}
\!\!\!
\mu(\tau)
\!\!\;\!\;\!
\int_{\mbox{}_{\scriptstyle \!\:\!
R\;\!+\;\!S \,<\, |\,x\,| \,<\, R \;\!+\;\! 2 \;\!S}}
\hspace{-2.440cm}
|\, \nabla u \,|^{\:\!p \;\!-\;\! 1}
\;\!\;\!
|\,\nabla \zeta_{\mbox{}_{\scriptstyle R, \:\!S}}(x)\,|
\;
dx \,d\tau,
\end{equation}
\mbox{} \vspace{-0.425cm} \\
\begin{equation}
\notag
\mbox{\footnotesize $H$}_{\mbox{}_{\scriptstyle \!\:\!1}}
\!\:\!(\mbox{\footnotesize $R$})
\;=\;\!
\int_{\mbox{}_{\scriptstyle 0}}^{\:\!T}
\!\!\;\!\:\!
\int_{\mbox{}_{\scriptstyle \!\:\!
R/2 \,<\, |\,x\,| \,<\, R}}
\hspace{-1.575cm}
|\, \mbox{\boldmath $g$}(\tau,u) \,|
\;
|\,\nabla \zeta_{\mbox{}_{\scriptstyle R, \:\!S}}(x)\,|
\;
dx \,d\tau,
\end{equation}
\mbox{} \vspace{-0.425cm} \\
\begin{equation}
\notag
\mbox{\footnotesize $H$}_{\mbox{}_{\scriptstyle \!\:\!2}}
\!\;\!(\mbox{\footnotesize $R$},\,\!\mbox{\footnotesize $S$})
\;=\;\!
\int_{\mbox{}_{\scriptstyle 0}}^{\:\!T}
\!\!\;\!\:\!
\int_{\mbox{}_{\scriptstyle \!\:\!
R\;\!+\;\!S \,<\, |\,x\,| \,<\, R \;\!+\;\! 2 \;\!S}}
\hspace{-2.440cm}
|\, \mbox{\boldmath $g$}(\tau,u) \,|
\;
|\,\nabla \zeta_{\mbox{}_{\scriptstyle R, \:\!S}}(x)\,|
\;
dx \,d\tau.
\end{equation}
\mbox{} \vspace{+0.050cm} \\
Recalling that
$ \:\!p \geq n $
(by hypothesis),
we observe that \\
\mbox{} \vspace{-0.600cm} \\
\begin{equation}
\notag
\mbox{\footnotesize $J$}_{\mbox{}_{\scriptstyle \!\:\!1}}
\!\;\!(\mbox{\footnotesize $R$})
\;\:\!\leq\;\;\!
\epsilon
\!\:\!
\int_{\mbox{}_{\scriptstyle 0}}^{\:\!T}
\!\!\!
\mu(\tau)
\!\!\;\!\;\!
\int_{\mbox{}_{\scriptstyle \!\:\!
R/2 \,<\, |\,x\,| \,<\, R }}
\hspace{-1.575cm}
|\,\nabla \zeta_{\mbox{}_{\scriptstyle R, \:\!S}}(x)\,|^{\:\!p}
\;
dx \,d\tau
\;\:\!+\;\:\!
\epsilon^{\scriptstyle \!\:\!
-\, \frac{\scriptstyle 1}{\scriptstyle \;\!p \;\!-\;\!1\:\!} }
\!\!\:\!
\int_{\mbox{}_{\scriptstyle 0}}^{\:\!T}
\!\!\!
\mu(\tau)
\!\!\;\!\;\!
\int_{\mbox{}_{\scriptstyle \!\:\!
R/2 \,<\, |\,x\,| \,<\, R }}
\hspace{-1.575cm}
|\,\nabla u \,|^{\:\!p}
\;\!\;\!
dx \,d\tau,
\end{equation}
\mbox{} \vspace{+0.025cm} \\
and similarly for
$ {\displaystyle
J_{\mbox{}_{\scriptstyle \!\;\!2}}
\!\;\!(\mbox{\footnotesize $R$},\,\!\mbox{\footnotesize $S$})
} $,
$ {\displaystyle
H_{\mbox{}_{\scriptstyle \!\:\!1}}
\!\:\!(\mbox{\footnotesize $R$})
} $
and
$ {\displaystyle
H_{\mbox{}_{\scriptstyle \!\;\!2}}
\!\;\!(\mbox{\footnotesize $R$},\,\!\mbox{\footnotesize $S$})
} $.
This gives,
letting
$ \:\!\mbox{\footnotesize $S$} \rightarrow \infty $, \\
\mbox{} \vspace{-0.100cm} \\
\mbox{} \hspace{+1.500cm}
$ {\displaystyle
\int_{\mbox{}_{\scriptstyle \!\:\!
|\,x\,| \,>\, R }}
\hspace{-0.750cm}
|\, u(x,t) \,|
\;\!\;\! dx
\;\leq\;\!
\int_{\mbox{}_{\scriptstyle \!\;\!
|\,x\,| \,>\, R/2}}
\hspace{-1.075cm}
|\, u_0(x) \,| \;\!\;\! dx
\;+\;
\mbox{\footnotesize $ {\displaystyle
\frac{\;\!2\;\!\mbox{\footnotesize $C$}}
     {\mbox{\footnotesize $R$}} }$}
\!\;\!
\int_{\mbox{}_{\scriptstyle 0}}^{\:\!T}
\!\!\!
\mbox{\footnotesize $F$}(\tau)
\!\!\;\!\;\!
\int_{\mbox{}_{\scriptstyle \!\:\!
|\,x\,| \,>\, R/2 }}
\hspace{-1.075cm}
|\, u(x,t) \,|^{\:\!\kappa \;\!+\;\! 1}
\,
dx \,d\tau
} $ \\
\mbox{} \vspace{+0.150cm} \\
\mbox{} \hspace{+4.600cm}
$ {\displaystyle
+\;\;
\epsilon^{\scriptstyle \!\:\!
-\, \frac{\scriptstyle 1}{\scriptstyle \;\!p \;\!-\;\!1\:\!} }
\!\!\:\!
\int_{\mbox{}_{\scriptstyle 0}}^{\:\!T}
\!\!\!
\mu(\tau)
\!\!\;\!\;\!
\int_{\mbox{}_{\scriptstyle \!\:\!
|\,x\,| \,>\, R/2 }}
\hspace{-1.075cm}
|\,\nabla u \,|^{\:\!p}
\;
dx \,d\tau
\,+\,
\mbox{\footnotesize $K$}_{\mbox{}_{\scriptstyle \!\:\!n}}
\!\;\!\;\!\epsilon
\;\!\;\!
\Bigl\{\;\! 1 + \!
\int_{\mbox{}_{\scriptstyle 0}}^{\:\!T}
\!\!\!
\mu(\tau)
\;\!\;\!
d\tau
\Bigr\}
} $ \\
\mbox{} \vspace{+0.150cm} \\
\mbox{} \hspace{+6.750cm}
$ {\displaystyle
+\;\;
\epsilon^{\scriptstyle \!\:\!
-\, \frac{\scriptstyle 1}{\scriptstyle \;\!p \;\!-\;\!1\:\!} }
\!\!\:\!
\int_{\mbox{}_{\scriptstyle 0}}^{\:\!T}
\!\!\;\!
\int_{\mbox{}_{\scriptstyle \!\:\!
|\,x\,| \,>\, R/2 }}
\hspace{-1.075cm}
|\, u(x,\tau) \,|
\;
dx \,d\tau
} $ \\
\mbox{} \vspace{+0.225cm} \\
for every
$ \;\!0 < t \leq \mbox{\footnotesize $T$} $,
where
$ {\displaystyle
\!\;\!\;\!
\mbox{\footnotesize $K$}_{\mbox{}_{\scriptstyle \!\;\!n}}
\!\;\!
} $
is some constant depending on $n$,
$\mbox{\footnotesize $C$}$ only
(and not on $ \mbox{\footnotesize $R$} $),
and where we have used (2.6)
and the assumption
$ \!\;\!\;\! p \geq n $.
Therefore,
by (1.3) and (2.3),
we can choose
$ \:\!\mbox{\footnotesize $R$} > 0 $
sufficiently large
(depending on $ \epsilon $,
$ \mbox{\footnotesize $T$} $)
such that \\
\mbox{} \vspace{-0.650cm} \\
\begin{equation}
\notag
\mbox{} \hspace{+1.000cm}
\int_{\mbox{}_{\scriptstyle \!\:\!
|\,x\,| \,>\, R }}
\hspace{-0.750cm}
|\, u(x,t) \,|
\;\!\;\! dx
\;\leq\;\;\!
\epsilon
\,+\;\!\;\!
\mbox{\footnotesize $K$}_{\mbox{}_{\scriptstyle \!\:\!n}}
\;\!\epsilon
\;\!\;\!
\Bigl\{\;\! 1 + \!
\int_{\mbox{}_{\scriptstyle 0}}^{\:\!T}
\!\!\!
\mu(\tau)
\;\!\;\!
d\tau
\Bigr\}
\;\;\qquad
\forall \;\,
0 \:\!<\:\! t \:\!\leq\:\!
\mbox{\footnotesize $T$}
\!\;\!.
\end{equation}
\mbox{} \vspace{-0.200cm} \\
Since $ \;\!\epsilon > 0 \:\!$
is arbitrary,
and the constant
$ \mbox{\footnotesize $K$}_{\mbox{}_{\scriptstyle \!\:\!n}} \!\:\!$
in the estimate above
is independent of $\:\!\epsilon$,
this gives
$ {\displaystyle
\;\!
u(\cdot,t)
\in
C^{0}(\,\![\,0, \;\!\mbox{\footnotesize $T$}_{\!\ast}),
\;\!L^{1}(\mathbb{R}^{n})\,\!)
} $,
as claimed in the first part of
\mbox{\small \sc Proposition 2.4}. \linebreak
\mbox{} \vspace{-0.750cm} \\

Finally,
to show the second part
(i.e., mass conservation),
we proceed in a similar way, \linebreak
but taking this time
$ \phi(x) = \zeta_{\mbox{}_{R}}\!\;\!(x) $
in (2.2),
where
$ \zeta_{\mbox{}_{R}}\!\;\!(\cdot) $
is the
cut-off function
considered in the proof
of \mbox{\footnotesize \sc Proposition 2.2}.
This completes the proof
of \mbox{\footnotesize \sc Proposition 2.4}.
\hfill $\Box$ \\
}
%
%
\nl
%
%
%
%
{\bf Remark 2.3.}
%
In a similar way,
in the remaining case $ \;\! p < n \:\!$
mass conservation
%
%
%
%
can be obtained
from (2.2) with
$ \phi(x) = \zeta_{\mbox{}_{R}}\!\;\!(x) $
provided that we have,
instead of (2.6),
the stronger condition \\
\mbox{} \vspace{-0.730cm} \\
\begin{equation}
\tag{2.9}
|\;\!\;\! \mbox{\boldmath $g$}(t,\mbox{u})\,|
\;\!\;\!\leq\,
C(\mbox{\small $T$})
\;\!\;\!
|\,\mbox{u}\,|
\qquad \;\;
\forall \;\;
|\,\mbox{u}\,| \:\!\ll\:\! 1,
\;\,
0 < t < \mbox{\small $T$}\!\:\!,
\end{equation}
\mbox{} \vspace{-0.215cm} \\
%
%
%
%
and that we have
$ {\displaystyle
\;\!
|\, \nabla u\:\!(\cdot,t) \,|
\,\mbox{\small $\in$}\,
L^{q}_{\mbox{\scriptsize loc}}
([\;\!0, \:\!\mbox{\small $T$}_{\!\ast}), \:\!
L^{q}(\mathbb{R}^{n})\,\!)
} $
for some
$ \;\! q \,\mbox{\small $\in$}\,[\;\!p - 1, p) $
satisfying
$ \;\!q \;\!\leq\:\! (\,\!p - 1)\;\!n\,\!/\,\!(\,\!n - 1) $.
For still other conditions,
see \cite{Guidolin2015}, Ch.\;2.

%
%
%
%

%
%
%
\newpage
\mbox{} \vspace{-1.000cm} \\
%
%
%
%
{\bf 3. \mbox{\boldmath $L^{1}\!\:\!$} contraction
and comparison properties} \\
\mbox{} \vspace{-0.600cm} \\

The results obtained in this section,
where we introduce a few extra assumptions
(see (3.1)\,-\,(3.4) below),
serve to establish the {\em uniqueness\/}
of solutions to (1.1), (1.2),
among other important properties
\cite{Kalashnikov1987, WuZhaoYinLi2001}.
Upon $ \!\:\!\mbox{\boldmath $f$} \!\:\!$
and
$ \mbox{\boldmath $g$} $,
it will be required
one of the following sets of conditions:
for every given
$ \:\!\mbox{\small $M$} > 0 $,
$ \:\!0 < \mbox{\small $T$} < \mbox{\small $T$}_{\!\ast} $,
one must have (1.6) and (1.7) satisfied,
that is, \\
\mbox{} \vspace{+0.000cm} \\
\mbox{} \hspace{+0.350cm}
$ {\displaystyle
|\,\mbox{\boldmath $f$}(x,t,\mbox{u})
\;\!-\:\!
\mbox{\boldmath $f$}(x,t,\mbox{v}) \,|
\;\leq\:
\mbox{\small $K$}_{\!\,\!f}\,\!
(\mbox{\small $M$}\!\;\!, \:\!\mbox{\small $T$})
\,
|\, \mbox{u} - \mbox{v} \,|^{\scriptstyle \:\!
1 \,-\, \frac{\scriptstyle 1}{\scriptstyle p} }
\quad \;\,
\forall \;\,
x \in \mathbb{R}^{n} \!\:\!,
\:
0 \:\!\leq\:\! t \leq\:\! \mbox{\small $T$}\!\;\!
} $,
\mbox{} \hfill (3.1) \\
\mbox{} \hspace{+9.900cm}
$ {\displaystyle
|\,\mbox{u}\,| \:\!\leq\:\! \mbox{\small $M$}\!\;\!,
\;
|\,\mbox{v}\,| \:\!\leq\:\! \mbox{\small $M$}\!\;\!
} $, \\
\mbox{} \vspace{-0.250cm} \\
\mbox{} \hspace{+1.000cm}
$ {\displaystyle
|\,\mbox{\boldmath $g$}(t,\mbox{u})
\;\!-\:\!
\mbox{\boldmath $g$}(t,\mbox{v}) \,|
\;\leq\,
\mbox{\small $K$}_{\!\,\!g}\,\!
(\mbox{\small $M$}\!\;\!, \:\!\mbox{\small $T$})
\;\!\;\!
|\, \mbox{u} - \mbox{v} \,|^{\scriptstyle \:\!
1 \,-\, \frac{\scriptstyle 1}{\scriptstyle p} }
\hspace{+1.650cm}
\forall \;\,
0 \:\!\leq\:\! t \leq\:\! \mbox{\small $T$}\!\;\!
} $,
\mbox{} \hfill (3.2) \\
\mbox{} \hspace{+9.900cm}
$ {\displaystyle
|\,\mbox{u}\,| \:\!\leq\:\! \mbox{\small $M$}\!\;\!,
\;
|\,\mbox{v}\,| \:\!\leq\:\! \mbox{\small $M$}\!\;\!
} $, \\
\mbox{} \vspace{-0.400cm} \\
or
the stronger assumptions \\
\mbox{} \vspace{-0.050cm} \\
\mbox{} \hspace{+1.150cm}
$ {\displaystyle
|\,\mbox{\boldmath $f$}_{\!\,\!\mbox{\scriptsize u}}
(x,t,\mbox{u})\,|
\;\leq\:
\mbox{\small $F$}_{\!\!\;\!\mbox{}_{\mbox{\scriptsize u}}}
(\mbox{\small $M$}\!\;\!, \:\!\mbox{\small $T$}) \,
|\, \mbox{u} \,|^{\:\!\kappa}
\qquad \;\,
\forall \;\,
x \in \mathbb{R}^{n} \!\:\!,
\;
0 \:\!\leq\:\! t \:\!\leq\:\! \mbox{\small $T$}\!\;\!,
\:
|\,\mbox{u}\,| \:\!\leq\:\! \mbox{\small $M$}\!\;\!
} $,
\mbox{} \hfill (3.3) \\
\mbox{} \vspace{-0.100cm} \\
\mbox{} \hspace{+1.400cm}
$ {\displaystyle
|\,\mbox{\boldmath $g$}_{\!\;\!\mbox{\scriptsize u}}
(t,\mbox{u})\,|
\;\leq\:
\mbox{\small $G$}_{\!\;\!\mbox{}_{\mbox{\scriptsize u}}}
\!\;\!(\mbox{\small $M$}\!\:\!, \:\!\mbox{\small $T$}) \,
|\, \mbox{u} \,|^{\:\!\gamma}
\hspace{+2.025cm}
\forall \;\,
0 \:\!\leq\:\! t \:\!\leq\:\! \mbox{\small $T$}\!\;\!,
\:
|\,\mbox{u}\,| \:\!\leq\:\! \mbox{\small $M$}\!\;\!
} $,
\mbox{} \hfill (3.4) \\
\mbox{} \vspace{-0.050cm} \\
with constants
$ \;\!
\mbox{\small $K$}_{\!\!\;\!f}\,\!
(\mbox{\small $M$}\!, \!\;\!\;\!\mbox{\small $T$}),
\;\!
\mbox{\small $K$}_{\!\!\;\!g}\,\!
(\mbox{\small $M$}\!, \!\;\!\;\!\mbox{\small $T$}),
\;\!
\mbox{\small $F$}_{\!\!\;\!\mbox{}_{\mbox{\scriptsize u}}}
\!\;\!(\mbox{\small $M$}\!, \!\;\!\;\!\mbox{\small $T$}),
\;\!
\mbox{\small $G$}_{\!\!\;\!\;\!\mbox{}_{\mbox{\scriptsize u}}}
\!\;\!(\mbox{\small $M$}\!, \!\;\!\;\!\mbox{\small $T$}) $
depending on
$ \:\!\mbox{\small $M$}\!\,\!\,\!, \;\!\mbox{\small $T$} \!\,\!$,
where
$ \mbox{\boldmath $f$}_{\!\,\!\mbox{\scriptsize u}}
\!\,\!=\:\! \partial \mbox{\boldmath $f$}
\!\;\!/\!\;\!\:\!\partial \:\!\mbox{u} $,
$ \mbox{\boldmath $g$}_{\!\;\!\mbox{\scriptsize u}}
\!\,\!=\:\! \partial \!\;\!\:\!\mbox{\boldmath $g$}
\!\;\!/\!\;\!\partial \:\!\mbox{u} $.
\!\;\!We note that
(3.3)\,-\,(3.4) are satisfied
in the prototype model
given by
$ \mbox{\boldmath $f$}(x,t,\mbox{u}) \!\;\!=\!\;\!\;\!
\mbox{\boldmath $b$}(x,t) \;\!|\;\!\mbox{u}\;\!|^{\:\!\kappa}
\!\;\!\;\!\mbox{u} $,
$ \;\!\mbox{\boldmath $g$}(t,\mbox{u}) \!\;\!=\!\;\!\;\!
\mbox{\boldmath $c$}(t) \;\!|\;\!\mbox{u}\;\!|^{\:\!\gamma}
\!\;\!\;\!\mbox{u} $. \\
\mbox{} \vspace{-0.850cm} \\

Again,
as in the previous section,
solutions to (1.1), (1.2)
are always meant
in the space
$ {\displaystyle
\,\!\,\!
C^{0}([\;\!0, \:\!\mbox{\small $T$}_{\!\ast}),
L^{1}_{\mbox{\scriptsize loc}}(\mathbb{R}^{n}))
} $
$ {\displaystyle
\cap\,
L^{p}_{\mbox{\scriptsize loc}}
([\;\!0, \:\!\mbox{\small $T$}_{\!\ast}),
W^{1\!\;\!,\:p}_{\mbox{\scriptsize loc}}(\mathbb{R}^{n}))
\!\;\!\cap
L^{\infty}_{\mbox{\scriptsize loc}}
([\;\!0, \:\!\mbox{\small $T$}_{\!\ast}), \:\! L^{1}(\mathbb{R}^{n})
\!\;\!\cap\! L^{\infty}(\mathbb{R}^{n}))
} $,
with its maximal
existence interval
given by
$ [\;\!0, \:\!\mbox{\small $T$}_{\!\ast} \!\;\!) $. \\
\nl
%
%
%
%
{\bf Proposition 3.1.}
\textit{%
Let
$ \;\!u(\cdot,t), \;\!v(\cdot,t) $,
$ 0 < t \leq \mbox{\small $T$}\!\;\!$,
be given
solutions of $\,\!(1.1a)$, $(1.2)$
corresponding to
initial states
$ {\displaystyle
\;\!
u_0, \:\!v_0 \in
L^{1}(\mathbb{R}^{n})
\cap
L^{\infty}(\mathbb{R}^{n})
} $,
respectively.
Then \\
}
\mbox{} \vspace{-0.675cm} \\
\begin{equation}
\tag{3.5}
\|\, u(\cdot,t) - v(\cdot,t) \,
\|_{\mbox{}_{\scriptstyle L^{1}(\mathbb{R}^{n})}}
\;\!\leq\;
\|\, u_0 - v_0 \;\!
\|_{\mbox{}_{\scriptstyle L^{1}(\mathbb{R}^{n})}}
\qquad
\forall \;\,
0 \:\!<\:\! t \:\!\leq\:\! \mbox{\small $T$}
\!\;\!,
\end{equation}
\mbox{} \vspace{-0.175cm} \\
\textit{%
provided that\,$:$
$($i\/$)$
$ p \geq n $,
and
$ \mbox{\boldmath $f$} \!$,
$ \mbox{\boldmath $g$} $
satisfy
$\!\;\!\;\!(3.1) \!\;\!$
and $\!\;\!\;\!(3.2)$ above,
or
$($when
$ \;\! 2 < p < n \!\;\!\;\!)\!\!:$
$($ii\/$)$
$ \;\!\kappa \geq 1 - 2/p $,
$ \;\!\gamma \geq 1 - 2/p $,
\:\!and
$ \:\!\mbox{\boldmath $f$} \!\:\!$,
$ \mbox{\boldmath $g$} \,\!$
satisfy
$\!\;\!\;\!(3.3) \!\;\!$
and
$\!\;\!\;\!(3.4) $,
respectively. \\
}
%
%
\nl
\mbox{} \vspace{-0.450cm} \\
%
%
%
{\small
{\bf Proof.}
Given $ \:\!h > 0 $,
$ \delta > 0 $,
$ \:\!\mbox{\footnotesize $R$} > 0 $,
let
$ \:\!\zeta_{\mbox{}_{R}} \!\in C^{1}(\mathbb{R}^{n}) \:\!$
be the cut-off function
considered in the proof of
\mbox{\footnotesize \sc Proposition 2.2}.
\!Let
$ \,\!u_{h} (\cdot,t) $,
$ v_{h} (\cdot,t) \,\!$
be the time Steklov regularizations
of $ u(\cdot,t) $, $ v(\cdot,t) $,
respectively.
Let
$ L_{\delta} \!\;\!\in C^{3}(\mathbb{R}^{n}) $
be defined
as in the proof of
\mbox{\footnotesize \sc Proposition 2.2},
and let
$ \:\!\theta(\cdot,t) \!\;\!:= u(\cdot,t) - v(\cdot,t) $,
$ \theta_{\!\;\!h} \!\;\!(\cdot,t) \!\;\!:=
u_{h} \!\;\!(\cdot,t) - v_{h} \!\;\!(\cdot,t) $.
\!Taking
$ {\displaystyle
\;\!
\phi(x) =
L_{\delta}^{\prime}(\theta_{\!\;\!h}\!\;\!(x,t))
\;\!
\zeta_{\mbox{}_{R}}\!\;\!(x)
} $
in the equations (2.2)
for
$ \:\!u_{h} \!\;\!(\cdot,t) $,
$ v_{h} \!\;\!(\cdot,t) $,
subtracting one from the other
and integrating the result
in the interval
$ (\,\!t_0, \:\!t\,\!) $,
where
$ \:\!0 < t_0 < t $,
we get,
letting
$ \:\! h \;\!\mbox{\scriptsize $\searrow$}\;\!0 \:\!$
and
$ \:\! t_0 \:\!\mbox{\scriptsize $\searrow$}\;\! 0 $,

\mbox{} \vspace{-0.750cm} \\
\mbox{} \hspace{+1.450cm}
$ {\displaystyle
\int_{\mbox{}_{\scriptstyle |\,x\,|\,<\,R}}
\hspace{-0.825cm}
L_{\delta}(\theta(x,t))
\,
\zeta_{\mbox{}_{R}}\!\;\!(x)
\, dx
\;+\:\!
\int_{\mbox{}_{\scriptstyle \!\;\!0}}^{\;\!t}
\!\!\:\!
\mu(\tau)
\!
\int_{\mbox{}_{\scriptstyle |\,x\,|\,<\,R}}
\hspace{-0.825cm}
L_{\delta}^{\prime\prime}(\theta)
\;\!\:\!
\langle \;\!\;\! \mbox{\boldmath $a$}(u,v),
\:\! \nabla \theta \;\!\;\!\rangle
\:
\zeta_{\mbox{}_{R}}\!\;\!(x)
\;\!\;\!
dx \,d\tau
} $ \\
\mbox{} \vspace{+0.300cm} \\
\mbox{} \hspace{+0.900cm}
$ {\displaystyle
\leq \,
\int_{\mbox{}_{\scriptstyle |\,x\,|\,<\,R}}
\hspace{-0.825cm}
L_{\delta}(\theta_{0}(x))
\,
\zeta_{\mbox{}_{R}}\!\;\!(x)
\, dx
\;+\:\!
\int_{\mbox{}_{\scriptstyle \!\;\!0}}^{\;\!t}
\!\!\:\!
\mu(\tau)
\!
\int_{\mbox{}_{\scriptstyle \!\:\!R/2 \,<\,|\,x\,|\,<\,R}}
\hspace{-1.625cm}
|\, L_{\delta}^{\prime}(\theta)\,|
\cdot
|\, \mbox{\boldmath $a$}(u,v) \,|
\cdot
|\, \nabla \zeta_{\mbox{}_{R}}\!\;\!(x) \,|
\; dx \,d\tau
\;\;\! +
} $ \\
\mbox{} \vspace{+0.275cm} \\
\mbox{} \hfill
$ {\displaystyle
\int_{\mbox{}_{\scriptstyle \!\;\!0}}^{\;\!t}
\!\!\:\!
\mu(\tau)
\!
\int_{\mbox{}_{\scriptstyle |\,x\,|\,<\,R}}
\hspace{-0.825cm}
L_{\delta}^{\prime\prime}(\theta)\,
\,
|\, [\,\mbox{\boldmath $\tilde{f}$}\,] \,|
\cdot
|\, \nabla \theta \,|
\;
\zeta_{\mbox{}_{R}}\!\;\!(x)
\; dx \,d\tau
\;\!\;\!+
\int_{\mbox{}_{\scriptstyle \!\;\!0}}^{\;\!t}
\!\!\:\!
\mu(\tau)
\!
\int_{\mbox{}_{\scriptstyle \!\,\! R/2\,<\,|\,x\,|\,<\,R}}
\hspace{-1.625cm}
|\, L_{\delta}^{\prime}(\theta) \,|
\cdot
|\, [\,\mbox{\boldmath $\tilde{f}$}\,] \,|
\cdot
|\, \nabla \zeta_{\mbox{}_{R}}\!\;\!(x) \,|
\; dx \,d\tau
} $ \\
\mbox{} \vspace{+0.350cm} \\
in view of (2.3),
where
$ \;\!\theta_{0} = u_{0} \!\;\!- v_{0} $,
$ \:\![\,\mbox{\boldmath $\tilde{f}$}\,] \equiv
[\,\mbox{\boldmath $\tilde{f}$}\,](x,\tau) =
\mbox{\boldmath $\tilde{f}$}(x,\tau,u(x,\tau)) -
\mbox{\boldmath $\tilde{f}$}(x,\tau,v(x,\tau)) $,
$ {\displaystyle
\mbox{\boldmath $\tilde{f}$} =
\mbox{\boldmath $f$} + \mbox{\boldmath $g$}
} $,
and
$ {\displaystyle
\;\!
\mbox{\boldmath $a$}(u,v) \:\!=\;\!
|\,\nabla u(x,\tau) \,|^{\:\!p - 2} \;\! \nabla u(x,\tau)
\,\!-\:\!
|\,\nabla v(x,\tau) \,|^{\:\!p - 2} \;\! \nabla v(x,\tau)
} $.
\!\;\!Noticing that \\
\mbox{} \vspace{+0.025cm} \\
\mbox{} \hspace{+1.250cm}
$ {\displaystyle
\langle \;\!\;\! \mbox{\boldmath $a$}(u,v),
\:\! \nabla \theta \;\!\;\!\rangle
\:=\;\:\!
\mbox{\footnotesize $ {\displaystyle \frac{1}{2} }$}
\;\!
\bigl(\:
|\, \nabla u\,|^{\:\!p - 2} +\;\!
|\, \nabla v\,|^{\:\!p - 2}
\;\!\bigr)
\,\:\!
|\, \nabla \theta \,|^{\:\!2}
\;\;\!+
} $ \\
\mbox{} \vspace{-0.300cm} \\
\mbox{} \hspace{+7.000cm}
$ {\displaystyle
+\;\,
\mbox{\footnotesize $ {\displaystyle
\frac{1}{2} }$}
\;\!
\bigl(\:
|\, \nabla u\,|^{\:\!p - 2} -\;\!
|\, \nabla v\,|^{\:\!p - 2}
\;\!\bigr)
\;\!
\bigl(\:
|\, \nabla u\,|^{\:\!2} -\;\!
|\, \nabla v\,|^{\:\!2}
\;\!\bigr)
} $ \\
\mbox{} \vspace{-0.300cm} \\
\mbox{} \hspace{+3.800cm}
$ {\displaystyle
\geq\;
\frac{1}{\;2^{\mbox{}^{\scriptstyle \:\!p - 1}}}
\:
|\, \nabla \theta \,|^{\:\!p}
} $ \\
\mbox{} \vspace{+0.125cm} \\
and
that
$ {\displaystyle
\;\!
|\, \mbox{\boldmath $a$}(u,v) \,|
\;\!\leq\,
|\,\nabla u \,|^{\:\!p - 1} +\;\!
|\,\nabla v \,|^{\:\!p - 1}
\!
} $,
we then have \\
\mbox{} \vspace{+0.025cm} \\
\mbox{} \hspace{+1.400cm}
$ {\displaystyle
\int_{\mbox{}_{\scriptstyle |\,x\,|\,<\,R}}
\hspace{-0.825cm}
L_{\delta}(\theta(x,t))
\,
\zeta_{\mbox{}_{R}}\!\;\!(x)
\, dx
\;+\;
\Bigl(\, 1 - \frac{2}{p} \,\Bigr)
\,
\frac{1}{\;2^{\mbox{}^{\scriptstyle \:\!p - 1}}}
\!
\int_{\mbox{}_{\scriptstyle \!\;\!0}}^{\;\!t}
\!\!\:\!
\mu(\tau)
\!
\int_{\mbox{}_{\scriptstyle |\,x\,|\,<\,R}}
\hspace{-0.825cm}
L_{\delta}^{\prime\prime}(\theta)
\;\!\:\!
|\, \nabla \theta \,|^{\:\!p}
\,
\zeta_{\mbox{}_{R}}\!\;\!(x)
\;\!\;\!
dx \,d\tau
} $ \\
\mbox{} \vspace{+0.200cm} \\
\mbox{} \hspace{+0.800cm}
$ {\displaystyle
\leq \;\!\;\!\;\!
\|\, u_{0} -\;\! v_{0} \;\!
\|_{\mbox{}_{\scriptstyle L^{1}(\mathbb{R}^{n})}}
\:\!+\!\;\!
\int_{\mbox{}_{\scriptstyle \!\;\!0}}^{\;\!t}
\!\!\:\!
\mu(\tau)
\!
\int_{\mbox{}_{\scriptstyle \!\:\!R/2 \,<\,|\,x\,|\,<\,R}}
\hspace{-1.625cm}
|\, L_{\delta}^{\prime}(\theta)\,|
\:
\bigl(\,
|\, \nabla u \,|^{\:\!p - 1} \!\;\!+\;\! |\,\nabla v\,|^{\:\!p - 1}
\:\!\bigr)
\,
|\, \nabla \zeta_{\mbox{}_{R}}\!\;\!(x) \,|
\;\,\!
dx \,d\tau
\;\;\! +
} $ \\
\mbox{} \vspace{+0.100cm} \\
\mbox{} \hfill
$ {\displaystyle
+\;\;\!\;\!
2 \!\,\!
\int_{\mbox{}_{\scriptstyle \!\;\!0}}^{\;\!t}
\!\!\:\!
\mu(\tau)^{\mbox{}^{\scriptstyle \!\!
-\, \frac{\scriptstyle 1}{\scriptstyle \;\!p\:\! - 1\:\!} }}
\!\!\!\:\!
\int_{\mbox{}_{\scriptstyle |\,x\,|\,<\,R}}
\hspace{-0.825cm}
L_{\delta}^{\prime\prime}(\theta)
\;\!\;\!
|\, [\,\mbox{\boldmath $f$}\,] \,|^{\mbox{}^{\scriptstyle
\frac{\scriptstyle p}{\scriptstyle \;\!p\:\! - 1\:\!} }}
\!\;\!
\zeta_{\mbox{}_{R}}\!\;\!(x)
\; dx \,d\tau
\,+
\int_{\mbox{}_{\scriptstyle \!\;\!0}}^{\;\!t}
\!\!\;\!
\int_{\mbox{}_{\scriptstyle \!\,\! R/2\,<\,|\,x\,|\,<\,R}}
\hspace{-1.625cm}
|\, L_{\delta}^{\prime}(\theta) \,|
\cdot
|\, [\,\mbox{\boldmath $f$}\,] \,|
\cdot
|\, \nabla \zeta_{\mbox{}_{R}}\!\;\!(x) \,|
\;
dx \,d\tau
} $ \\
\mbox{} \vspace{+0.100cm} \\
\mbox{} \hfill
$ {\displaystyle
+\;\;\!\!\;\!\;\!
2 \!\,\!
\int_{\mbox{}_{\scriptstyle \!\;\!0}}^{\;\!t}
\!\!\:\!
\mu(\tau)^{\mbox{}^{\scriptstyle \!\!
-\, \frac{\scriptstyle 1}{\scriptstyle \;\!p\:\! - 1\:\!} }}
\!\!\!\:\!
\int_{\mbox{}_{\scriptstyle |\,x\,|\,<\,R}}
\hspace{-0.825cm}
L_{\delta}^{\prime\prime}(\theta)
\;\!\;\!
|\, [\,\mbox{\boldmath $g$}\,] \,|^{\mbox{}^{\scriptstyle
\frac{\scriptstyle p}{\scriptstyle \;\!p\:\! - 1\:\!} }}
\!\;\!
\zeta_{\mbox{}_{R}}\!\;\!(x)
\; dx \,d\tau
\,+\!\;\!
\int_{\mbox{}_{\scriptstyle \!\;\!0}}^{\;\!t}
\!\!\;\!
\int_{\mbox{}_{\scriptstyle \!\,\! R/2\,<\,|\,x\,|\,<\,R}}
\hspace{-1.625cm}
|\, L_{\delta}^{\prime}(\theta) \,|
\cdot
|\, [\,\mbox{\boldmath $g$}\,] \,|
\cdot
|\, \nabla \zeta_{\mbox{}_{R}}\!\;\!(x) \,|
\;
dx \,d\tau
} $, \\
\mbox{} \vspace{-0.350cm} \\
\mbox{} \hfill (3.6) \\
\mbox{} \vspace{-0.250cm} \\
where,
as before,
$ \:\![\,\mbox{\boldmath $f$}\,] \equiv
[\,\mbox{\boldmath $f$}\,](x,\tau) =
\mbox{\boldmath $f$}(x,\tau,u(x,\tau)) -
\mbox{\boldmath $f$}(x,\tau,v(x,\tau)) $,
$ \:\![\,\mbox{\boldmath $g$}\,] \equiv
[\,\mbox{\boldmath $g$}\,](x,\tau) =
\mbox{\boldmath $g$}(\tau,u(x,\tau)) -
\mbox{\boldmath $g$}(\tau,v(x,\tau)) $.
If $ \;\! p \geq n $,
we may proceed
as in the proof of
{\footnotesize \sc Proposition 2.4}
(using that
$ |\, L_{\delta}^{\prime}(\vartheta) \,|
\leq 1 $
for any $\:\!\vartheta \in \mathbb{R}$),
letting
$ \;\!\delta \rightarrow 0 \;\!$
and then
$ \:\!\mbox{\footnotesize $R$} \rightarrow \infty \:\!$
to obtain,
given
$ \;\! \epsilon > 0 \:\!$
arbitrary$\:\!$: \\
\mbox{} \vspace{-0.625cm} \\
\begin{equation}
\notag
\|\, \theta(\cdot,t) \,
\|_{\mbox{}_{\scriptstyle L^{1}(\mathbb{R}^{n})}}
\;\!\;\!\leq\;\;\!
\|\, u_{0} -\;\! v_{0} \;\!
\|_{\mbox{}_{\scriptstyle L^{1}(\mathbb{R}^{n})}}
\;\!+\:
\mbox{\footnotesize $K$}_{\mbox{}_{\scriptstyle \!\:\!n}}
\;\!\epsilon
\;\!\;\!
\Bigl\{\;\! 1 \:\!+
\!
\int_{\mbox{}_{\scriptstyle 0}}^{\:\!T}
\!\!\!
\mu(\tau)
\;\!\;\!
d\tau
\,\!\Bigr\}
\end{equation}
\mbox{} \vspace{-0.150cm} \\
for each $ \;\!0 < t \leq \mbox{\footnotesize $T$} \!\;\!$,
because of
(1.3), (2.3) and (3.1), (3.2) above,
where
$ \!\;\!\;\!\mbox{\small $K$}_{\mbox{}_{\scriptstyle \!\:\!n}} \! > 0 \:\!$
is some appropriate constant
depending on the dimension $\,\!\,\!n\,\!\,\!$ but not on $\,\!\,\!\epsilon$.
Since
this holds for any $ \epsilon > 0 $,
(3.5) is obtained
in the case $ \:\! p \geq n $,
as claimed. \\
\mbox{} \vspace{-0.775cm} \\
%

%
%

When
$ \;\! 2 < p < n $,
we assume (3.3), (3.4)
with
$ \kappa \geq 0 $,
$ \gamma \geq 0 $
satisfying
$ \kappa \geq 1 - 2/p $
and
$ \gamma \geq 1 - 2/p $,
proceeding
instead
as follows.
Because
$ {\displaystyle
\;\!
|\, L_{\delta}^{\prime}(\theta) \,|
\!\;\!\;\!\leq
K \;\!|\;\!\theta\;\!|\;\!/\;\!\delta
\;\!
} $
for all
$ \theta \in \mathbb{R} $,
$ \delta > 0 $
(and some constant
$ \!\;\!K $
independent of
\mbox{$ \!\;\!\;\!\theta, \delta $}),
we obtain,
letting
$ \mbox{\footnotesize $R$} \rightarrow \infty \:\!$
in (3.6):

\mbox{} \vspace{-0.950cm} \\
\begin{equation}
\notag
\int_{\mbox{}_{\scriptstyle \mathbb{R}^{n}}}
\!\!\!\!\;\!
L_{\delta}(\theta(x,t))
\;
dx
\;\,\!\leq\;\:\!
\|\, u_0 - \!\;\!\;\!v_0 \;\!
\|_{\mbox{}_{\scriptstyle L^{1}(\mathbb{R}^{n})}}
+\;
2 \!\:\!
\int_{\mbox{}_{\scriptstyle \!\;\!0}}^{\;\!t}
\!\!\:\!
\mu(\tau)^{\mbox{}^{\scriptstyle \!\!
-\, \frac{\scriptstyle 1}{\scriptstyle \;\!p\:\! - 1\:\!} }}
\!\!\!\:\!
\int_{\mbox{}_{\scriptstyle \mathbb{R}^{n}}}
\!\!\!\!\;\!
L_{\delta}^{\prime\prime}(\theta)
\;
|\, [\,\mbox{\boldmath $f$}\,](x,\tau) \,
|^{\mbox{}^{\scriptstyle
\frac{\scriptstyle p}{\scriptstyle \;\!p\:\! - 1\:\!} }}
\;\!
dx \,d\tau
\end{equation}
\mbox{} \vspace{-0.200cm} \\
\mbox{} \hspace{+6.500cm}
$ {\displaystyle
+\;\,
2 \!\:\!
\int_{\mbox{}_{\scriptstyle \!\;\!0}}^{\;\!t}
\!\!\:\!
\mu(\tau)^{\mbox{}^{\scriptstyle \!\!
-\, \frac{\scriptstyle 1}{\scriptstyle \;\!p\:\! - 1\:\!} }}
\!\!\!\:\!
\int_{\mbox{}_{\scriptstyle \mathbb{R}^{n}}}
\!\!\!\!\;\!
L_{\delta}^{\prime\prime}(\theta)
\;
|\, [\,\mbox{\boldmath $g$}\,](x,\tau) \,
|^{\mbox{}^{\scriptstyle
\frac{\scriptstyle p}{\scriptstyle \;\!p\:\! - 1\:\!} }}
\;\!
dx \,d\tau
} $ \\
\mbox{} \hfill (3.7) \\
\mbox{} \vspace{-0.450cm} \\
by (1.3) and (2.3).
Now,
because of (3.3) and (3.4),
we have \\
\mbox{} \vspace{-0.500cm} \\
\begin{equation}
\notag
|\, [\,\mbox{\boldmath $f$}\,](x,\tau) \,|
\;\leq\;\!\;\!
\sqrt{\:\!n\,}
\,
\mbox{\footnotesize $F$}_{\mbox{}_{\!\mbox{\scriptsize u}}}
\!\;\!
(\mbox{\footnotesize $M$}\!\:\!, \:\!\mbox{\footnotesize $T$})
\;\!\;\!
\bigl(\:
|\,u(x,\tau)\,|^{\:\!\kappa}
\!\;\!+\;\!
|\,v(x,\tau)\,|^{\:\!\kappa}
\;\!\bigr)
\;\!\;\!
|\,\theta(x,\tau)\,|,
\end{equation}
\mbox{} \vspace{-0.800cm} \\
\begin{equation}
\notag
|\, [\,\mbox{\boldmath $g$}\,](x,\tau) \,|
\;\leq\;\!\;\!
\sqrt{\:\!n\,}
\,
\mbox{\footnotesize $G$}_{\mbox{}_{\!\;\!\mbox{\scriptsize u}}}
\!\;\!
(\mbox{\footnotesize $M$}\!\:\!, \:\!\mbox{\footnotesize $T$})
\;\!\;\!
\bigl(\:
|\,u(x,\tau)\,|^{\:\!\gamma}
\!\;\!+\;\!
|\,v(x,\tau)\,|^{\:\!\gamma}
\;\!\bigr)
\;\!\;\!
|\,\theta(x,\tau)\,|
\end{equation}
\mbox{} \vspace{-0.150cm} \\
for all $ \;\!x \in \mathbb{R}^{n} \!\;\!$,
$ 0 < \tau \leq \mbox{\footnotesize $T$} \!\:\!$,
\;\!where
$ {\displaystyle
\;\!
\mbox{\footnotesize $M$}
=\;\!
\sup\,\bigl\{\;\!
\|\, u(\cdot,\tau) \,
\|_{L^{\infty}(\mathbb{R}^{n})}
\!\;\!,
\,
\|\, v(\cdot,\tau) \,
\|_{L^{\infty}(\mathbb{R}^{n})}
\!\!\;\!:\,
0 < \tau \leq \mbox{\footnotesize $T$}
\!\;\!\;\!\bigr\}
} $,
so that \\
\mbox{} \vspace{-0.750cm} \\
\begin{equation}
\notag
|\, [\,\mbox{\boldmath $f$}\,](x,\tau) \,
|^{\mbox{}^{\scriptstyle
\frac{\scriptstyle p}{\scriptstyle \;\!p\;\! - 1\:\!} }}
\:\!\leq\;\!\;\!
\mbox{\footnotesize $K$}\!\;\!
(\mbox{\footnotesize $M$}\!\:\!,
\:\!\mbox{\footnotesize $T$}\!\:\!, \;\!p, \:\!n)
\;\!\;\!
\biggl\{\,
|\,u(x,\tau)\,
|^{\mbox{}^{\scriptstyle \!
\frac{\scriptstyle \;\!\kappa \:\!p \;\!+\:\!1\:\!}
     {\scriptstyle \;\!p\;\! - 1\:\!} }}
\!\!+\;\!\;\!
|\,v(x,\tau)\,
|^{\mbox{}^{\scriptstyle \!
\frac{\scriptstyle \;\!\kappa \;\!p \;\!+\:\!1\:\!}
     {\scriptstyle \;\!p\;\! - 1\:\!} }}
\biggr\}
\;\!\;\!
|\,\theta(x,\tau)\,|
\end{equation}
\mbox{} \vspace{-0.450cm} \\
and \\
\mbox{} \vspace{-1.050cm} \\
\begin{equation}
\notag
|\, [\,\mbox{\boldmath $g$}\,](x,\tau) \,
|^{\mbox{}^{\scriptstyle
\frac{\scriptstyle p}{\scriptstyle \;\!p\;\! - 1\:\!} }}
\:\!\leq\;\!\;\!
\mbox{\footnotesize $K$}\!\;\!
(\mbox{\footnotesize $M$}\!\:\!,
\:\!\mbox{\footnotesize $T$}\!\:\!, \;\!p, \:\!n)
\;\!\;\!
\biggl\{\,
|\,u(x,\tau)\,
|^{\mbox{}^{\scriptstyle \!
\frac{\scriptstyle \;\!\gamma \:\!p \;\!+\:\!1\:\!}
     {\scriptstyle \;\!p\;\! - 1\:\!} }}
\!\!+\;\!\;\!
|\,v(x,\tau)\,
|^{\mbox{}^{\scriptstyle \!
\frac{\scriptstyle \;\!\gamma \;\!p \;\!+\:\!1\:\!}
     {\scriptstyle \;\!p\;\! - 1\:\!} }}
\biggr\}
\;\!\;\!
|\,\theta(x,\tau)\,|
\end{equation}
\mbox{} \vspace{-0.200cm} \\
for all $\:\!(x,\tau) \:\!$
concerned,
\:\!where
$ {\displaystyle
\;\!
\mbox{\footnotesize $K$}\!\;\!
(\mbox{\footnotesize $M$}\!\:\!,
\:\!\mbox{\footnotesize $T$}\!\:\!, \;\!p, \:\!n)
} $
is some constant
that does not depend on $\;\!\delta $.
Hence,
letting
$ \;\!\delta \rightarrow 0 \;\!$
in (3.7),
we obtain \\
\mbox{} \vspace{-0.650cm} \\
\begin{equation}
\notag
\|\, \theta(\cdot,t) \,
\|_{\mbox{}_{\scriptstyle L^{1}(\mathbb{R}^{n})}}
\;\!\leq\;
\|\, u_{0} -\:\! v_{0} \;\!
\|_{\mbox{}_{\scriptstyle L^{1}(\mathbb{R}^{n})}}
\end{equation}
\mbox{} \vspace{-0.250cm} \\
by Lebesgue's dominated convergence,
since
$ {\displaystyle
\;\!(\kappa\:\!p + 1)/(p - 1)
\geq 1
} $,
$ {\displaystyle
(\gamma\:\!p + 1)/(p - 1)
\geq 1
} $.
\!\;\!This shows (3.5)
in case ({\em ii\/}),
so that the proof of
\mbox{\footnotesize \sc Proposition 3.1}
is now complete.
\mbox{} \hfill $\Box$ \\
}
%
%
\mbox{} \vspace{-0.525cm} \\

Actually, under the same assumptions
of \mbox{\small \sc Proposition 3.1},
a lot more is true, \linebreak
as shown by the next two results
(cf. \mbox{\small \sc Propositions 3.2}
and \mbox{\small \sc 3.3} below): \\
\nl
%
%
%
%
{\bf Proposition 3.2.}
\textit{%
Let
$ \;\!u(\cdot,t), \;\!v(\cdot,t) $,
$ 0 < t \leq \mbox{\small $T$}\!\;\!$,
be given
solutions of $\,\!(1.1a)$, $(1.2)$
corresponding to
initial states
$ {\displaystyle
\;\!
u_0, \:\!v_0 \in
L^{1}(\mathbb{R}^{n})
\cap
L^{\infty}(\mathbb{R}^{n})
} $,
respectively.
Then \\
}
\mbox{} \vspace{-0.675cm} \\
\begin{equation}
\tag{3.8}
\bigl\|\;\!\;\! \bigl(\:\! u(\cdot,t) - v(\cdot,t)\:\!\bigr)_{\mbox{}_{\!\!+}}
\bigr\|_{\mbox{}_{\scriptstyle \!\;\!L^{1}(\mathbb{R}^{n})}}
\,\!\leq\;
\bigl\|\;\!\;\! \bigl(\:\! u_0 - v_0 \:\!\bigr)_{\mbox{}_{\!\!+}}
\bigr\|_{\mbox{}_{\scriptstyle \!\;\!L^{1}(\mathbb{R}^{n})}}
\qquad
\forall \;\,
0 \:\!<\:\! t \:\!\leq\:\! \mbox{\small $T$}
\end{equation}
{\em and} \\
\mbox{} \vspace{-0.925cm} \\
\begin{equation}
\tag{3.9}
\bigl\|\;\!\;\! \bigl(\:\! u(\cdot,t) - v(\cdot,t)\:\!\bigr)_{\mbox{}_{\!\!-}}
\bigr\|_{\mbox{}_{\scriptstyle \!\;\!L^{1}(\mathbb{R}^{n})}}
\,\!\leq\;
\bigl\|\;\!\;\! \bigl(\:\! u_0 - v_0 \:\!\bigr)_{\mbox{}_{\!\!-}}
\bigr\|_{\mbox{}_{\scriptstyle \!\;\!L^{1}(\mathbb{R}^{n})}}
\qquad
\forall \;\,
0 \:\!<\:\! t \:\!\leq\:\! \mbox{\small $T$}
\!\;\!,
\end{equation}
\mbox{} \vspace{-0.125cm} \\
\textit{%
provided that\,$:$
$($i\/$)$
$ p \geq n $,
and
$ \mbox{\boldmath $f$} \!$,
$ \mbox{\boldmath $g$} $
satisfy
$\!\;\!\;\!(3.1) \!\;\!$
and $\!\;\!\;\!(3.2)$ above,
or
$($when
$ \;\! 2 < p < n \!\;\!\;\!)\!\!:$
$($ii\/$)$
$ \;\!\kappa \geq 1 - 2/p $,
$ \;\!\gamma \geq 1 - 2/p $,
\:\!and
$ \:\!\mbox{\boldmath $f$} \!\:\!$,
$ \mbox{\boldmath $g$} \,\!$
satisfy
$\!\;\!\;\!(3.3) \!\;\!$
and
$\!\;\!\;\!(3.4) $,
respectively. \\
}
%
%
(\:\!Here, as usual,
$ \:\!\theta_{\mbox{}_{\!+}} \!$
and
$ \:\!\theta_{\mbox{}_{\!-}} \!$
stand for the positive and negative real parts,
respectively,
of a given number $ \theta \in \mathbb{R} $,
that is\:\!:
$ \theta_{\mbox{}_{\!+}} \!=\:\! (\;\!|\;\!\theta\;\!| + \theta\;\!)/2 $,
and
$ \;\!\theta_{\mbox{}_{\!-}} \!=\:\! (\;\!|\;\!\theta\;\!| - \theta\;\!)/2 $.) \\
\nl
\mbox{} \vspace{-0.450cm} \\
%
%
%
{\small
{\bf Proof.}
The following argument is adapted
from the proof of \mbox{\footnotesize \sc Proposition 3.1}
and \cite{DiehlFabrisZingano2014, SchutzZiebellZinganoZingano2013}:
taking $ H \in C^{2}(\mathbb{R}) $
such that
$ H^{\prime}(\xi) \geq 0 $ for all $\xi \in \mathbb{R}$,
$ H(\xi) = 0 $ $\,\forall \;\!\;\!\xi \leq 0 $,
$ H(\xi) = 1 $ $\forall \;\xi \geq 1 $,
and given $ \delta > 0 $ (arbitrary),
let $ H_{\mbox{}_{\scriptstyle \!\:\!\delta}} \!\in C^{2}(\mathbb{R}) $
be defined by
$ H_{\mbox{}_{\scriptstyle \!\:\!\delta}}(\xi) \!\:\!:= H(\xi/\delta) $.
Also,
given \linebreak
$ \:\!h > 0 $,
$ \:\!\mbox{\footnotesize $R$} > 0 $,
let
$ \:\!\zeta_{\mbox{}_{R}} \!\in C^{1}(\mathbb{R}^{n}) \:\!$
be the cut-off function
used in the proof of
\mbox{\footnotesize \sc Proposition 2.2}.
Letting
$ u_{\mbox{}_{\scriptstyle h}} \!(\cdot,t) $,
$ v_{\mbox{}_{\scriptstyle h}} \!(\cdot,t) $
denote the Steklov regularizations
of $ u(\cdot,t) $, $ v(\cdot,t) $,
respectively, \linebreak
and setting
$ \:\!\theta(\cdot,t) \!\:\!:= u(\cdot,t) - v(\cdot,t) $,
$ \theta_{\mbox{}_{\scriptstyle \!\;\! h}} \!(\cdot,t) \!\:\!:=
u_{\mbox{}_{\scriptstyle \!\;\! h}} \!(\cdot,t) -
v_{\mbox{}_{\scriptstyle \!\;\! h}} \!(\cdot,t) $,
we may proceed as follows.
Taking
$ {\displaystyle
\;\!
\phi(x) =
H_{\mbox{}_{\scriptstyle \!\:\!\delta}}
\!\;\!(\theta_{\mbox{}_{\scriptstyle \!\;\!h}}\!\:\!(x,t))
\;\!
\zeta_{\mbox{}_{R}}\!\;\!(x)
} $
in the equations (2.2)
for
$ \:\!u_{\mbox{}_{\scriptstyle \!\;\! h}} \!(\cdot,t) $,
$ v_{\mbox{}_{\scriptstyle \!\;\! h}} \!(\cdot,t) $,
subtracting one from the other
and integrating the result
in the interval
$ (\,\!t_0, \:\!t\,\!) $,
where
$ \:\!0 < t_0 < t $,
we get,
letting
$ \:\! h \;\!\mbox{\scriptsize $\searrow$}\;\!0 \:\!$
and
$ \:\! t_0 \:\!\mbox{\scriptsize $\searrow$}\;\! 0 $, \\
\mbox{} \vspace{+0.050cm} \\
\mbox{} \hspace{+1.450cm}
$ {\displaystyle
\int_{\mbox{}_{\scriptstyle |\,x\,|\,<\,R}}
\hspace{-0.825cm}
G_{\mbox{}_{\scriptstyle \!\:\!\delta}}
\!\;\!(\:\!\theta(x,t))
\;\!\;\!
\zeta_{\mbox{}_{R}}\!\;\!(x)
\, dx
\;+\:\!
\int_{\mbox{}_{\scriptstyle \!\;\!0}}^{\;\!t}
\!\!\:\!
\mu(\tau)
\!
\int_{\mbox{}_{\scriptstyle |\,x\,|\,<\,R}}
\hspace{-0.825cm}
H_{\!\:\!\delta}^{\,\prime}(\theta)
\;\!\:\!
\langle \;\!\;\! \mbox{\boldmath $a$}(u,v),
\:\! \nabla \theta \;\!\;\!\rangle
\:
\zeta_{\mbox{}_{R}}\!\;\!(x)
\;\!\;\!
dx \,d\tau
} $ \\
\mbox{} \vspace{+0.300cm} \\
\mbox{} \hspace{+0.900cm}
$ {\displaystyle
\leq \,
\int_{\mbox{}_{\scriptstyle |\,x\,|\,<\,R}}
\hspace{-0.825cm}
G_{\mbox{}_{\scriptstyle \!\:\!\delta}}
\!\;\!(\theta_{0}(x))
\;\!\;\!
\zeta_{\mbox{}_{R}}\!\;\!(x)
\, dx
\;+\:\!
\int_{\mbox{}_{\scriptstyle \!\;\!0}}^{\;\!t}
\!\!\:\!
\mu(\tau)
\!
\int_{\mbox{}_{\scriptstyle \!\:\!R/2 \,<\,|\,x\,|\,<\,R}}
\hspace{-1.625cm}
|\, H_{\mbox{}_{\scriptstyle \!\:\!\delta}}(\theta)\,|
\cdot
|\, \mbox{\boldmath $a$}(u,v) \,|
\cdot
|\, \nabla \zeta_{\mbox{}_{R}}\!\;\!(x) \,|
\; dx \,d\tau
\;\;\! +
} $ \\
\mbox{} \vspace{+0.275cm} \\
\mbox{} \hfill
$ {\displaystyle
\int_{\mbox{}_{\scriptstyle \!\;\!0}}^{\;\!t}
\!\!\:\!
\mu(\tau)
\!
\int_{\mbox{}_{\scriptstyle |\,x\,|\,<\,R}}
\hspace{-0.825cm}
H_{\!\:\!\delta}^{\,\prime}(\theta)\,
\,
|\, [\,\mbox{\boldmath $\tilde{f}$}\,] \,|
\cdot
|\, \nabla \theta \,|
\;
\zeta_{\mbox{}_{R}}\!\;\!(x)
\; dx \,d\tau
\;\!\;\!+
\int_{\mbox{}_{\scriptstyle \!\;\!0}}^{\;\!t}
\!\!\:\!
\mu(\tau)
\!
\int_{\mbox{}_{\scriptstyle \!\,\! R/2\,<\,|\,x\,|\,<\,R}}
\hspace{-1.625cm}
|\, H_{\mbox{}_{\scriptstyle \!\:\!\delta}}(\theta) \,|
\cdot
|\, [\,\mbox{\boldmath $\tilde{f}$}\,] \,|
\cdot
|\, \nabla \zeta_{\mbox{}_{R}}\!\;\!(x) \,|
\; dx \,d\tau
\!\;\!,
} $ \\
\mbox{} \vspace{+0.350cm} \\
where
$ \;\!G_{\mbox{}_{\scriptstyle \!\:\!\delta}}\!\;\!(\eta)
= \int_{0}^{\:\!\eta}
H_{\mbox{}_{\scriptstyle \!\:\!\delta}}\!\;\!(\xi)
\, d\xi $,
$ \;\! [\,\mbox{\boldmath $\tilde{f}$}\,] \equiv
[\,\mbox{\boldmath $\tilde{f}$}\,](x,\tau) =
\mbox{\boldmath $\tilde{f}$}(x,\tau,u(x,\tau)) -
\mbox{\boldmath $\tilde{f}$}(x,\tau,v(x,\tau)) $,
$ {\displaystyle
\mbox{\boldmath $\tilde{f}$} \!:=\!\;\!
\mbox{\boldmath $f$} \!+ \mbox{\boldmath $g$}
} $,
and
$ {\displaystyle
\;\!
\mbox{\boldmath $a$}(u,v) \:\!=
|\,\nabla u(x,\tau) \,|^{\:\!p - 2} \;\! \nabla u(x,\tau)
\,\!-
|\,\nabla v(x,\tau) \,|^{\:\!p - 2} \;\! \nabla v(x,\tau)
} $,
as before.
\!\!\;\!From this point,
we repeat the steps in the proof of
\mbox{\footnotesize \sc Proposition 3.1},
using now that
$ \;\!G_{\mbox{}_{\scriptstyle \!\;\!\delta}}\!(\theta)
\rightarrow \theta_{\mbox{}_{\!\:\!+}} \!\;\!$
as $ \delta \rightarrow 0 \!\;\!\;\!$:
in case ({\em i\/}),
we let
$ \delta \rightarrow 0 $ and
$ \mbox{\footnotesize $R$} \rightarrow \infty $
to obtain (3.8),
and
in case ({\em ii\/})
we reverse the order,
letting this time
$ \mbox{\footnotesize $R$} \rightarrow \infty $
and then
$ \delta \rightarrow 0 $
to arrive at (3.8),
as claimed.

The proof of (3.9) follows exactly the same lines,
except that this we take
$ H \in C^{2}(\mathbb{R}) $
satisfying\:\!:
$ H^{\prime}(\xi) \leq 0 $ for all $\xi \in \mathbb{R}$,
$ H(\xi) = 1 $ $\,\forall \;\!\;\!\xi \leq -1 $,
and
$ H(\xi) = 0 $ $\forall \;\xi \geq 0 $.
\mbox{} \hfill $\Box$ \\
}
%
\mbox{} \vspace{-0.475cm} \\

A direct consequence of (3.8) (or of (3.9))
is the following comparison principle.
\nl
%
%
%
%
{\bf Proposition 3.3.}
\textit{%
Let
$ \;\!u(\cdot,t), \;\!v(\cdot,t) $,
$ 0 < t \leq \mbox{\small $T$}\!\;\!$,
be given
solutions of $\,\!(1.1a)$, $(1.2)$
corresponding to
initial states
$ {\displaystyle
\;\!
u_0, \:\!v_0 \in
L^{1}(\mathbb{R}^{n})
\cap
L^{\infty}(\mathbb{R}^{n})
} $,
respectively.
Then \\
}
\mbox{} \vspace{-0.550cm} \\
\begin{equation}
\tag{3.10}
u_0(\;\!\cdot\;\!)
\,\leq\,
v_0(\;\!\cdot\;\!)
\;\;\;\;\;
\Longrightarrow
\;\;\;\;\;
u(\cdot,t)
\,\leq\,
v(\cdot,t)
\quad \;\;
\forall \;\,
0 \:\!<\:\! t \:\!\leq\:\! \mbox{\small $T$}
\!\;\!,
\end{equation}
\mbox{} \vspace{-0.175cm} \\
\textit{%
provided that\,$:$
$($i\/$)$
$ p \geq n $,
and
$ \mbox{\boldmath $f$} \!$,
$ \mbox{\boldmath $g$} $
satisfy
$\!\;\!\;\!(3.1) \!\;\!$
and $\!\;\!\;\!(3.2)$ above,
or
$($when
$ \;\! 2 < p < n \!\;\!\;\!)\!\!:$
$($ii\/$)$
$ \;\!\kappa \geq 1 - 2/p $,
$ \;\!\gamma \geq 1 - 2/p $,
\:\!and
$ \:\!\mbox{\boldmath $f$} \!\:\!$,
$ \mbox{\boldmath $g$} \,\!$
satisfy
$\!\;\!\;\!(3.3) \!\;\!$
and
$\!\;\!\;\!(3.4) $,
respectively. \\
}
%
%

\newpage

%
%

\mbox{} \vspace{-1.200cm} \\

{\bf Acknowledgements.}
This work was partly supported by
{\small CNP}q (\mbox{\small M}inistry of \mbox{\small S}cience
and \mbox{\small T}echnology, \mbox{\small B}razil),
\mbox{\small G}rant \mbox{\small \#\,154037/2011-7}
and by {\small CAPES}
(\mbox{\small M}inistry of \mbox{\small E}ducation,
\mbox{\small B}razil),
\mbox{\small G}rant
\mbox{\small \#\,1212003/2013}.
The authors also express
their gratitude to Paulo R. Zingano
(\mbox{\small UFRGS}, \mbox{\small B}razil)
for some helpful suggestions
and discussions. \linebreak
\mbox{} \vspace{-0.750cm} \\
%
%
%
%
%
%
%

%
%

%
\nl
\nl
\nl
{\small
\begin{minipage}[t]{10.00cm}
\mbox{\normalsize \textsc{Jocemar de Quadros Chagas}} \\
Departamento de Matem\'atica e Estat\'\i stica \\
Universidade Estadual de Ponta Grossa \\
Ponta Grossa, PR 84030-900, Brazil \\
E-mail: {\sf jocemarchagas@uepg.br} \\
\end{minipage}
\nl
\mbox{} \vspace{-0.450cm} \\
\nl
\begin{minipage}[t]{10.00cm}
\mbox{\normalsize \textsc{Patr\'\i cia Lisandra Guidolin}} \\
Instituto Federal de Educa\c c\~ao, Ci\^encia e Tecnologia \\
Farroupilha, RS 95180-000, Brazil \\
E-mail: {\sf patricia.guidolin@farroupilha.ifrs.edu.br} \\
\end{minipage}
\nl
\mbox{} \vspace{-0.450cm} \\
\nl
\begin{minipage}[t]{10.00cm}
\mbox{\normalsize \textsc{Jana\'\i na Pires Zingano}} \\
Departamento de Matem\'atica Pura e Aplicada \\
Universidade Federal do Rio Grande do Sul \\
Porto Alegre, RS 91509-900, Brazil \\
E-mail: {\sf jzingano@mat.ufrgs.br} \\
\end{minipage}
}
%
%

\end{document}